\documentclass[leqno,12pt]{article}

\usepackage{theorem}
\usepackage{amssymb}
\usepackage{amsmath}
\usepackage{amsfonts}
\usepackage{times}

\usepackage[mathcal]{euler}
\usepackage{euler}

\def \RR {\mathbb R}

\def \EE {\mathbb E}
\def \ZZ {\mathbb Z}
\def \TT {\mathbb T}
\def \CC {\mathbb C}
\def \PP {\mathbb P}

\def \vphi {\varphi}

\def \cD {\mathcal D}

\def \cS {\mathcal S}

\def \cG {\mathcal G}
\def \cH {\mathcal H}

\newtheorem{theorem}{Theorem}[section]
\newtheorem{lemma}[theorem]{Lemma}

\newtheorem{proposition}[theorem]{Proposition}
\newtheorem{corollary}[theorem]{Corollary}
 {\theorembodyfont{\rmfamily}}
\newtheorem{definition}[theorem]{Definition}
\newtheorem{remark}[theorem]{Remark}

\def\myffrac#1#2 in #3{\raise 2.6pt\hbox{$#3 #1$}\mkern-1.5mu\raise 0.8pt\hbox{$
#3/$}\mkern-1.1mu\lower 1.5pt\hbox{$#3 #2$}}

\begin{document}

\title{Poincar\'e Inequalities and Moment Maps}
\author{Bo'az Klartag\thanks{{School of Mathematical Sciences, Tel-Aviv University, Tel Aviv
69978, Israel. Supported in part by the Israel Science Foundation
and by a Marie Curie Reintegration Grant from the Commission of the
European Communities. Email: klartagb@tau.ac.il} }}
\date{}
\maketitle

\abstract{We propose a new method for obtaining Poincar\'e-type
inequalities on arbitrary convex bodies in $\RR^n$. Our technique
involves a dual version of Bochner's formula and a certain moment
map, and it also applies to some non-convex sets. In particular, we
generalize the central limit theorem for convex bodies to a  class
of non-convex domains, including the unit balls of $\ell_p$-spaces
in $\RR^n$ for $0 < p < 1$.}

\section{Introduction}

An important observation that goes back to Sudakov \cite{sudakov}
and to Diaconis and Freedman \cite{DF} is that approximately
gaussian marginals are intimately related to {\it thin shell
inequalities}. That is, let $X$ be a random vector in $\RR^n$ with
mean zero and identity covariance, where the dimension $n$ is
assumed very high. Suppose that $X$ satisfies a  thin shell
inequality, of the form
\begin{equation}  \EE \left( \frac{|X|^2}{n} -1 \right)^2 \ll
1. \label{eq_153}
\end{equation} It then follows that there are plenty of vectors
$\theta \in \RR^n$ for which the scalar product $\langle X, \theta
\rangle$ is approximately a gaussian random variable. See  von
Weizs\"acker \cite{vW}, Bobkov \cite{bobkov}, Anttila, Ball and
Perissinaki \cite{abp} or \cite{K_clt, K_euro} for further
explanations, and Eldan and Klartag \cite{EK} for connections to the
hyperplane conjecture.

\medskip In this paper, Poincar\'e-type inequalities refer to
inequalities in which the variance of a function is bounded in terms
of an integral of a quadratic form involving the gradient of the
function. One of the methods used to prove a thin shell bound such as
(\ref{eq_153}) goes through such Poincar\'e-type inequalities in
high-dimensional spaces. This approach was pursued in \cite{ptrf},
where the Bochner formula was applied to study optimal thin shell
bounds and Poincar\'e-type inequalities for the uniform measure on
high-dimensional convex bodies. The technique in \cite{ptrf} and in
the related work by Barthe and Cordero-Erausquin \cite{bce} relied
very much on symmetries of the probability distribution under
consideration. The method seemed quite irrelevant for arbitrary convex
bodies, possessing no symmetries. The following twist is proposed
here: Introduce additional symmetries by considering a certain
transportation of measure from a space of twice or thrice the
dimension.  The plan is to apply Bochner's formula in this higher
dimensional space, and deduce a Poincar\'e-type inequality for the
original measure.

\medskip We proceed by demonstrating the Poincar\'e-type
inequalities that are obtained in the simplest case, perhaps, in
which the convex set we investigate is $\RR^n_{+}$, the orthant of
all $x \in \RR^n$ with positive coordinates. A function $\vphi:
\RR^n_{+} \rightarrow (-\infty, \infty]$ is called $p$-convex, for
$0 < p \leq 1$, if the function
$$  (x_1,\ldots,x_n)  \mapsto
\vphi \left(x_1^{1/p},\ldots,x_n^{1/p} \right) $$ is convex on
$\RR^n_+$. For instance $\vphi(x) = \sum_{i=1}^n \sqrt{x_i}$ is
$p$-convex for any $p \leq 1/2$.

 \begin{theorem} Let $n \geq 1, k > 1$ be integers. Suppose that
 $\mu$ is a Borel measure on $\RR^n_{+}$
with density $\exp(-\vphi)$, where $\vphi: \RR^n_{+} \rightarrow
(-\infty, \infty]$ is $p$-convex for $p=1/k$. Assume  that $f:
\RR^n_+ \rightarrow \RR$ is a $\mu$-integrable, locally Lipschitz
function with $\int f d \mu = 0$. Then,
\begin{equation}
 \int_{\RR^n_{+}} f^2 d \mu \leq \frac{k^2}{k-1} \sum_{i=1}^n
\int_{\RR^n_{+}} x_i^{2} \left| \partial^i f(x) \right|^2 d \mu(x).
\label{eq_1015} \end{equation} \label{almost_gap} Here, $\partial^i
f = \partial f / \partial x_i$ stands for the derivative of $f$ with
respect to the $i^{th}$ variable.
\end{theorem}
We emphasize that the function $f$ in Theorem \ref{almost_gap} is not
assumed to satisfy any boundary conditions. Compare, for example,
to the Hardy-type inequalities in
 Matskewich and Sobolevskii \cite{MS}. We say that a subset $K
 \subset \RR^n_+$ is $p$-convex for $0 < p \leq 1$, if
 $$ \left \{ (x_1^p,\ldots,x_n^p)  \ ; \ (x_1,\ldots,x_n) \in K
 \right \} $$
 is a convex set. In other words, $K$ is $p$-convex when the function
 that equals $0$ on $K$ and equals $+\infty$ outside $K$ is
 $p$-convex. Observe that the intersection of $p$-convex sets is again a $p$-convex
 set. Dilations centered at the origin preserve $p$-convexity. For $p \neq 1$, translations do not necessarily preserve
 $p$-convexity, but $p$-convexity is preserved by translations conjugated with the map $x \mapsto (x_1^p,\ldots,x_n^p)$.
  From Theorem \ref{almost_gap} we immediately deduce:

 \begin{corollary} Let $n \geq 1, \ell > 1$ be integers, and assume that $K
 \subset \RR^n_+$ is a $(1/\ell)$-convex set with a non-empty interior. Then, for
 any locally Lipschitz, integrable function $f: K \rightarrow \RR$
 with $\int_K f = 0$,
$$ \int_K f^2 dx \leq \frac{\ell^2}{\ell-1} \sum_{i=1}^n
\int_{K} x_i^{2} \left| \partial^i f(x) \right|^2 d x.
$$ \label{cor_510}
 \end{corollary}

 For $x, y  \in \RR^n_{+}$ we write $x \leq y$ when $x_i \leq y_i$
for $i=1,\ldots,n$. A function $\vphi: \RR^n_{+} \rightarrow
(-\infty, \infty]$ is {\it increasing} when
$$ x \leq y \quad \quad \Longrightarrow \quad\quad
\vphi(x) \leq \vphi(y) \quad\quad\quad\quad\quad\quad (\text{for} \
x, y \in \RR^n_{+}). $$ It is simple to see that when $f$ is
increasing and $p$-convex, it is also $q$-convex for any $0 < q <
p$. A convex function is obviously $1$-convex. A function $\vphi:
\RR^n \rightarrow (-\infty, \infty]$ is unconditional if
$$ \vphi(x_1,\ldots,x_n) = \vphi(|x_1|, \ldots, |x_n|) \quad \quad \quad \quad (x \in \RR^n). $$
Observe that when $\vphi$ is an unconditional, convex function on
$\RR^n$, the restriction $\vphi|_{\RR^n_{+}}$ is necessarily
increasing and $p$-convex for any $0 < p \leq 1$. Thus Corollary
\ref{cor_510} recovers the Poincar\'e-type inequalities from
\cite{ptrf}: Quite unexpectedly, the unconditionality is used only
to infer that when $\vphi|_{\RR^n_+}$ is $1$-convex,  it is also
$(1/2)$-convex. Theorem \ref{almost_gap} may be generalized to
measures on $\RR^n$ whose density is unconditional, as follows:

 \begin{theorem} Let $\mu$ be a probability measure on $\RR^n$
with density $\exp(-\vphi)$, where $\vphi: \RR^n \rightarrow
(-\infty, \infty]$ is unconditional, and $\vphi|_{\RR^n_+}$ is
increasing and $1/k$-convex for an integer $k > 1$. Denote
$$ V_i = \int_{\RR^n} x_i^2 d \mu(x) \quad \quad \quad (i=1,\ldots,n). $$
Then, for any $\mu$-integrable, locally Lipschitz function $f: \RR^n
\rightarrow \RR$ with $\int f d \mu = 0$,
\begin{equation}
 \int_{\RR^n} f^2 d \mu \leq
\int_{\RR^n} \sum_{i=1}^n \left( \frac{k^2}{k-1} x_i^{2} + V_{i}
\right) \left|
\partial^i f(x) \right|^2 d \mu(x). \label{eq_0037}
\end{equation} \label{thm_easy}
Furthermore, when the function $f$ is unconditional, we may
eliminate the $V_i$'s on the right-hand side of (\ref{eq_0037}).
\end{theorem}

For $0 < p < 1$, denote by $\mu_p$  the uniform probability measure
on the non-convex set
$$ B_p^{n} = \left \{ x \in \RR^n \  ; \ \sum_{i=1}^n |x_i|^p \leq 1 \right \}. $$
Theorem \ref{thm_easy} applies for the measure $\mu_p$, with $k =
\lceil 1/p \rceil$.  Substituting  $f(x) = |x|^2 - \int |y|^2 d\mu(y)$
into Theorem \ref{thm_easy} yields thin shell bounds, which may be
used to infer the existence of approximately gaussian marginals. Further discussion of the central limit theorem for {\it
fractionally-convex bodies}, such as those in Theorem
\ref{thm_easy}, is deferred to a future work. Once Theorem
\ref{almost_gap} and Corollary \ref{cor_510} are formulated, one is
tempted to try and find a more direct proof of these inequalities.
In Section \ref{aftermath} we discuss such a direct argument, based
on the Brascamp-Lieb inequality \cite{BL}, and obtain
generalizations of Theorem \ref{almost_gap} and Theorem
\ref{thm_easy} in which $k
> 1$ is not necessarily an integer. Similarly, $\ell > 1$ does
not have to be an integer in Corollary \ref{cor_510}.

\medskip Next, suppose $K \subset \RR^n$ is a convex body, i.e., a bounded, open
convex set. We turn to the details of the Poincar\'e-type
inequalities that are obtained for $K$. Recall that a function on
$\RR^n$ is log-concave if it takes the form $\exp(-H)$ for a convex
function $H:\RR^n \rightarrow (-\infty, \infty]$. A Borel measure on
$\RR^n$ is log-concave if its density is log-concave, and in
particular, the uniform probability measure on an open, convex set
is log-concave.
 We say that a smooth, convex function $\psi: \RR^n \rightarrow \RR$ induces
a ``log-concave transportation to $K$'' if the following two
conditions hold:
\begin{enumerate}
\item[(a)] The function $\rho_{\psi}(x) = \det \nabla^2 \psi(x)$ is positive and log-concave
on $\RR^n$, where $\nabla^2 \psi$ is the Hessian of $\psi$.
\item[(b)] We have $\nabla \psi(\RR^n)= K$, where $\nabla
\psi(\RR^n) = \{ \nabla \psi(x) ; x \in \RR^n \}$.
\end{enumerate}
Observe that the map $x \mapsto \nabla \psi(x)$ pushes forward the
measure whose density is $\rho_\psi$, to the uniform measure on the
convex body $K$. For a given convex body $K \subset \RR^n$, there
are plenty of convex functions $\psi$ that induce a log-concave
transportation to $K$. In fact, for any log-concave function $\rho$
on $\RR^n$ whose integral equals the volume of $K$, there exists a
convex function $\psi$ which satisfies (a) and (b) with $\rho_{\psi}
= \rho$. This follows from the general theory of optimal
transportation of measure (e.g., Villani \cite{villani}). For
indices $i,j,k=1,\ldots,n$ we abbreviate
$$ \psi_i = \frac{\partial \psi}{\partial x_i}, \quad \psi_{ij} = \frac{\partial^2 \psi}{\partial x_i \partial x_j},
\quad \psi_{ijk} = \frac{\partial^3 \psi}{\partial x_i \partial x_j
\partial x_k}. $$
We also write $\left( \psi^{ij} \right)_{i,j=1,\ldots,n}$ for the
inverse matrix to the Hessian matrix $\nabla^2 \psi = \left(
\psi_{ij} \right)_{i,j=1,\ldots,n}$. The Legendre transform of
$\psi$ is the function $\psi^*: K \rightarrow \RR$ defined via
$$ \psi^*(x) = \sup_{y \in \RR^n} \left[ \langle x, y \rangle - \psi(y)
\right]. $$ Then $\nabla \psi^*$ is the inverse map to $\nabla
\psi$. With any $x \in K$ we associate the  quadratic form
$Q^{*}_{\psi, x}$ on $\RR^n$ defined by $$ Q^{*}_{\psi, x}(V) =
\sum_{i,j,k,\ell, m,p=1}^n V^i V^j \psi^{\ell m} \psi_{jkm} \psi^{k
p} \psi_{i \ell p} $$
 where $V =
(V^1,\ldots,V^n) \in \RR^n$ and where the functions $\psi_{ij},
\psi^{\ell m}, \psi_{jkm}$ etc. are evaluated at the point $\nabla
\psi^*(x)$. For $x \in K$ and $U \in \RR^n$, set
$$ Q_{\psi,x}(U) = \sup \left \{ 4 \left( \sum_{i,j=1}^n \psi_{ij} U^i V^j \right)^2 \ ; \
 V \in \RR^n, \ Q^*_{\psi, x}(V) \leq 1 \right \}, $$
 where $\psi_{ij}$ is evaluated at the point $\nabla \psi^*(x)$.
It could occur that $Q_{\psi, x}(U)$ is finite only for $U$ in a
certain subspace $E \subset \RR^n$. Note that $Q_{\psi, x}$ is a
quadratic form on that subspace $E$.

 \medskip There is
 one technical assumption that we must make. In Section
 \ref{sec3} we define the notion of  {\it regularity at infinity} of the function $\psi$,
 and throughout the analysis below we conveniently assume the $\psi$ is indeed regular at
 infinity.  This assumption seems to
 hold in the examples that we consider.  In the case where $K \subset \RR^n$ is a simple rational polytope,
regularity at infinity was investigated by
 Abreu \cite{abr_orbi}, who explained that it holds under fairly
 mild assumptions.

\begin{theorem} Let $K \subset \RR^n$ be a convex body. Suppose that
$\psi: \RR^n \rightarrow \RR$ induces a log-concave transportation
to $K$. Assume further that $\psi$ is regular at infinity. Then, for
any Lipschitz function $f: K \rightarrow \RR$,
$$ \int_K f = 0 \quad \quad \Rightarrow \quad \quad \int_K f^2 \leq
\int_K Q_{\psi,x}( \nabla f(x) ) dx. $$ \label{thm_2253}
\end{theorem}

 In order to apply
Theorem \ref{thm_2253}. one needs to select a function $\psi$ which induces a
log-concave transportation to $K$. Unfortunately, we are currently
unaware of a general method for constructing a ``reasonable''
function $\psi$ that satisfies (a) and (b), with good control over
derivatives up to order three. In simple cases, such as when $K
\subset \RR^n$ is the cube or the simplex, Theorem \ref{thm_2253}
does yield meaningful inequalities. See Section \ref{sec4} for a
detailed analysis of the case of the simplex. In particular, Theorem
\ref{thm_513} below provides somewhat unusual Poincar\'e-type
inequalities for a class of distributions on the regular simplex. We
present the proof of Theorem \ref{almost_gap} in Section \ref{sec2},
before dealing with the more general Theorem \ref{thm_2253} in
Section \ref{sec3}. In Section \ref{sec5} we prove Theorem
\ref{thm_easy}. Throughout this paper, by a {\it smooth function} we
mean a $C^{\infty}$-smooth one.

\medskip {\it{Acknowledgements.}} Thanks to Semyon Alesker, Franck Barthe, Dmitry
Faifman, Uri Grupel, Greg Kuperberg, Emanuel Milman, Yaron Ostrover,
Leonid Polterovich, Yanir Rubinstein and Mikhail Sodin for
interesting related discussions.

\section{Non-Linear Measure Projection}
\label{sec2}

In this section we prove Theorem \ref{almost_gap}. The analysis in
this section is also intended to serve as a preparation for Section
\ref{sec3}.
 Let $n,k \geq 1$ be positive integers, fixed throughout
this section. Denote $m = nk$. We use
$$ z = (z_1,\ldots,z_n) \in (\RR^k)^n = \RR^{kn} $$
as coordinates in $\RR^{kn}$, where $z_1,\ldots,z_n$ are
$k$-dimensional vectors. Consider the map $\pi: \RR^{m} \rightarrow
\overline{\RR^n_{+}}$ defined by
$$ \pi(z) = (|z_1|^k,\ldots,|z_n|^k) \quad \quad \quad \quad (z_1,\ldots,z_n) \in (\RR^k)^n. $$
Here, $\overline{\RR^n_{+}}$ is the closure of $\RR^n_+$ in $\RR^n$,
and $|z_i|$ stands for the standard Euclidean norm of $z_i \in
\RR^k$. The continuous map $\pi$ is proper, meaning that
$\pi^{-1}(K)$ is compact whenever $K \subset \overline{\RR_{+}^n}$
is compact. Let $S^{k-1} = \{ y \in \RR^k ; |y| = 1 \}$ denote the
unit sphere in $\RR^k$, and more generally, let $S^{k-1}(R) = \{ y
\in \RR^k ; |y| = R \}$. We write $\sigma_R$ for the uniform
probability measure on the sphere $S^{k-1}(R)$. With any $x \in
\RR^n_{+}$ we associate the cartesian product of spheres,
$$ \pi^{-1}(x) := S^{k-1}(x_1^{1/k}) \times S^{k-1}(x_2^{1/k}) \times \ldots \times
S^{k-1}(x_n^{1/k}) \subseteq (\RR^k)^n = \RR^m. $$ We denote by
$\sigma_x$ the uniform probability measure on $\pi^{-1}(x)$, that
is, the direct product of the uniform probability measures on the
spheres $S^{k-1}(x_j^{1/k})$ for $j=1,\ldots,n$.

\medskip We view the map $\pi$ as a kind of  {\it moment map}. The case $k=2$ fits very well
with the standard terminology, as in this case $\pi$ is related to
the moment map associated with the symplectic action of the group
$\left( SO(2) \right)^n$ on $(\RR^2)^n$ (see, e.g., Cannas da Silva
\cite{c}).  In the following lemma we verify that indeed the uniform
measure on $\RR^{m}$ is pushed forward to the uniform measure on
$\RR^n_{+}$ via the map $\pi$, up to a normalizing coefficient. We
write $Vol_k$ for the standard $k$-dimensional volume measure.

\begin{lemma} For any integrable function $f: \RR^n_{+} \rightarrow \RR$,
\begin{equation}
 \int_{\RR^{m}} f(\pi(z)) d Vol_m(z) = \omega_{n,k} \int_{\RR^n_{+}}
f(x) d Vol_n(x) \label{eq_545} \end{equation} where  $\omega_{n,k} =
\left( \pi^{k/2} / \Gamma(k/2 + 1) \right)^n$ is the $n^{th}$ power
of the volume of the $k$-dimensional unit ball. Furthermore, for any
Borel set $A \subseteq \RR^m$,
\begin{equation} Vol_m(A) = \omega_{n,k} \int_{\RR^n_{+}} \sigma_x(A)
d Vol_n(x). \label{eq_546} \end{equation}  \label{lem_2311}
\end{lemma}

\emph{Proof:} Integrating in polar coordinates for each $z_j \in \RR^k \ (j=1,\ldots,n)$, we find that
$$
 \int_{\RR^{m}} f(|z_1|^k,\ldots,|z_n|^k) dz_1 \ldots d z_n =
\omega_k^n \int_{\RR^n_{+}} f(x_1^k,\ldots,x_n^k) \left(
\prod_{j=1}^n
 x_j^{k-1} \right) dx_1 \ldots dx_n, $$where $\omega_k = k \pi^{k/2}
/ \Gamma(k/2 + 1)$ is the surface area of the unit sphere in
$\RR^k$. Applying the change of variables $(t_1,\ldots,t_n) =
(x_1^k,\ldots,x_n^k)$ we obtain
$$ \int_{\RR^n_{+}} f(x_1^k,\ldots,x_n^k) \left(
\prod_{j=1}^n x_j^{k-1} \right) dx_1 \ldots dx_n = k^{-n}
\int_{\RR^n_{+}} f(t_1,\ldots,t_n) dt_1 \ldots dt_n $$ and
(\ref{eq_545}) follows. The relation (\ref{eq_546}) is proven in a
similar fashion. \hfill $\square$

\medskip  Suppose $\nu$ is a Borel
measure on $\RR^{m}$. For a function $f \in L^2(\nu)$ we define
\begin{equation} \| f \|_{H^{-1}(\nu)} = \sup \left \{ \int_{\RR^{m}} f g d \nu
\, ; \, \int_{\RR^{m}} |\nabla g|^2 d \nu \leq 1 \right \},
\label{eq_1636} \end{equation} where the supremum runs over all
smooth functions $g: \RR^m \rightarrow \RR$ that belong to
$L^2(\nu)$. Note that $\| f \|_{H^{-1}(\nu)} = + \infty$ when $\int
f d \nu \neq 0$. The square of the $H^{-1}(\nu)$-norm is
sub-additive in $\nu$, as will be proven next:

\begin{lemma} Suppose $\nu$ is a  Borel
measure on $\RR^{m}$ that takes the form \begin{equation} \nu =
\int_{\Omega} \nu_{\alpha} d \lambda(\alpha) \label{eq_1535}
\end{equation}
 for  Borel
measures $\{ \nu_{\alpha} \}_{\alpha \in \Omega}$ on $\RR^m$ and a
measure $\lambda$ on $\Omega$. Then, for any $f \in L^2(\nu)$,
$$ \| f \|_{H^{-1}(\nu)}^2 \leq \int_{\Omega} \| f
\|_{H^{-1}(\nu_\alpha)}^2 d \lambda(\alpha). $$ \label{lem_1009}
\end{lemma}

 \emph{Proof:} Let $g$ be a  smooth function on $\RR^m$ which belongs to $L^2(\nu)$.
Since $f, g \in L^2(\nu_{\alpha})$ for $\lambda$-almost any $\alpha
\in \Omega$, then
 $$    \left| \int_{\RR^m} f g d \nu_\alpha \right| \leq  \| f
 \|_{H^{-1}(\nu_\alpha)} \sqrt{ \int_{\RR^m} \left| \nabla g
 \right|^2 d \nu_{\alpha}} $$
for $\lambda$-almost any $\alpha \in \Omega$.
 From (\ref{eq_1535}) and the Cauchy-Schwartz inequality,
 \begin{align*}
  \left| \int_{\RR^m} f g d \nu \right| & \leq \int_{\Omega}   \| f
 \|_{H^{-1}(\nu_\alpha)} \left( \int_{\RR^m} \left| \nabla g
 \right|^2 d \nu_{\alpha} \right)^{1/2}  d \lambda(\alpha) \\ & \leq \sqrt{ \int_{\RR^m} \| f
 \|_{H^{-1}(\nu_\alpha)}^2 d \lambda(\alpha)  } \cdot \sqrt{ \int_{\RR^m} |\nabla g|^2 d \nu} . \end{align*}
\hfill $\square$

\medskip
Recall that we use $(z_1,\ldots,z_n) \in (\RR^k)^n$ as coordinates
in $\RR^m = \RR^{kn}$. Let us furthermore denote $z_\ell =
(z_\ell^1,\ldots,z_\ell^k) \in \RR^k$, for any $\ell=1,\ldots,n$.

\begin{lemma} Assume $k \geq 2$. Let $x \in
\RR^n_{+}$. Let $1 \leq \ell \leq n, 1 \leq j \leq k$, and denote
$f(z) = z_{\ell}^j$ for $z \in \RR^m$. Then,
$$ \left \| f \right \|_{H^{-1}(\sigma_x)} \leq \frac{x_{\ell}^{2/k}}{\sqrt{k(k-1)}}. $$
\label{lem_2234}
\end{lemma}

\emph{Proof:} We claim that for any smooth function $h: \RR^k
\rightarrow \RR$ and $\theta \in S^{k-1}$,
\begin{equation}
 \int_{S^{k-1}}  \left \langle y,  \theta \right \rangle h(y) d \sigma_1(y) \leq \sqrt{\frac{1}{k (
 k-1)}} \cdot
 \sqrt{ \int_{S^{k-1}} |\nabla
h|^2 d \sigma_1 }. \label{eq_2201}
\end{equation}
Indeed, (\ref{eq_2201}) simply expresses the standard fact that $y
\mapsto \sqrt{k} (y \cdot \theta)$ is a normalized eigenfunction of
the Laplace-Beltrami operator on $S^{k-1}$, corresponding to the
eigenvalue $k-1$ (see, e.g., M\"uller \cite{muller}). By scaling, we
see that for any $R > 0$ and $\theta \in S^{k-1}$,
\begin{equation}
 \int_{S^{k-1}(R)} \left \langle y,  \theta \right \rangle h(y) d \sigma_R(y) \leq \frac{R^2}{\sqrt{k (
 k-1)}} \cdot
 \sqrt{ \int_{S^{k-1}(R)} |\nabla
h|^2 d \sigma_{R} }. \label{eq_2220}
\end{equation}
According to (\ref{eq_2220}), for any fixed $z_1,\ldots,z_{\ell-1},
z_{\ell+1},\ldots,z_n \in \RR^k$ and a smooth function $g: \RR^m
\rightarrow \RR$,
$$ \int_{S^{k-1}(R_{\ell})} z_{\ell}^j g(z_1,\ldots,z_{n}) d
\sigma_{R_\ell}(z_\ell) \leq \frac{x_{\ell}^{2/k}}{\sqrt{k(k-1)}}
\sqrt{ \int_{S^{k-1}(R_{\ell})} |\nabla g(z)|^2 d
\sigma_{R_{\ell}}(z_\ell)}, $$ where $R_{\ell} = x_{\ell}^{1/k}$.
Recall that the probability measure $\sigma_{x}$ is a product
measure, and that $\sigma_{R_{\ell}}$ is the $\ell^{th}$ factor in
this product. Integrating with respect to the remaining variables
$z_1,\ldots,z_{\ell-1}, z_{\ell+1},\ldots,z_n$, and using the
Cauchy-Schwartz inequality, we obtain
$$ \int_{\pi^{-1}(x)} z_{\ell}^j
g(z) d \sigma_{x}(z) \leq \frac{x_{\ell}^{2/k}}{\sqrt{k(k-1)}}
\sqrt{ \int_{\pi^{-1}(x)} |\nabla g(z)|^2 d \sigma_x(z)}. $$ The lemma
follows from the definition of the $H^{-1}(\sigma_x)$-norm. \hfill
$\square$

\medskip The following lemma is one of the reasons for considering the higher-dimensional
space $\RR^m$, rather than working in the original space
$\RR^n_{+}$. The extra dimensions translate to ``extra symmetries'',
which substitute for the explicit symmetries assumed in
\cite[Corollary 5]{ptrf} and in Barthe and Cordero-Erausquin
\cite[Section 3]{bce}. This effect actually seems more prominent in
Section \ref{sec3}.

\begin{lemma} Assume $k \geq 2$, let $1 \leq \ell \leq n, 1 \leq j \leq k$ and let $x \in
\RR^n_{+}$. Suppose that $f: \RR^n_{+} \rightarrow \RR$ is
differentiable at $x$. Denote $g(z) = f(\pi(z))$ for $z \in \RR^m$.
Then,
$$
 \left \| \frac{\partial g}{\partial z_\ell^j} \right
\|_{H^{-1}(\sigma_x)} \leq \sqrt{ \frac{k}{k-1} } \cdot x_{\ell}
 \left|
\partial^\ell f (x) \right|.
$$\label{lem_2141}
\end{lemma}

\emph{Proof:} Note that for $z \in \pi^{-1}(x)$,
$$ \frac{\partial g}{\partial z_\ell^j} (z_1,\ldots,z_n) =
k |z_\ell|^{k-2} z_\ell^j  \cdot \partial^\ell
f(|z_1|^k,\ldots,|z_n|^k) = \left( k x_{\ell}^{(k-2) / k}
\partial^\ell f (x_1,\ldots,x_n) \right) z_\ell^j.
$$
That is, the function $\left. \partial g \right/ \partial z_\ell^j$
is proportional to the linear function $z \mapsto z_\ell^j$ on the
support of $\sigma_x$, and the proportion coefficient is exactly $ k
x_{\ell}^{(k-2) / k}
\partial^\ell f (x_1,\ldots,x_n) $. According to Lemma
\ref{lem_2234},
\begin{align*}
 \left \| \frac{\partial g}{\partial z_\ell^j}  \right
\|_{H^{-1}(\sigma_x)}  & =  k x_{\ell}^{(k-2) / k}  \left|
\partial^\ell f (x_1,\ldots,x_n) \right| \cdot \left \| z_{\ell}^j \right
\|_{H^{-1}(\sigma_x)} \\ & \leq k x_{\ell}^{(k-2) / k} \left|
\partial^\ell f (x_1,\ldots,x_n) \right| \cdot \frac{x_{\ell}^{2/k}}{\sqrt{k(k-1)}}. \end{align*}
\hfill $\square$

\medskip Suppose $\Omega \subset \RR^m$ is a bounded, open set. We say that a
smooth function $u: \Omega \rightarrow \RR$ is {\it smooth up to the
boundary} if all of its derivatives of all orders are bounded in
$\Omega$. Note that when $u$ is smooth up to the boundary,  the
boundary values of $u$ and its derivatives are well-defined on
$\partial \Omega$, by continuity. For $R
> 1$  denote
$$ \Omega_{R} = \left \{ (z_1,\ldots,z_n) \in (\RR^k)^n \ ; \ R^{-1} < |z_i| < R \
\ \text{for} \ \ i=1,\ldots, n \right \}. $$ We denote by
$\partial_{reg} \Omega_R$ the regular part of the boundary $\partial
\Omega_R$. That is,  $$
\partial_{reg} \Omega_R = \left( \bigcup_{i=1}^n A_i^- \right) \cup
\left( \bigcup_{i=1}^n A_i^+ \right) $$ where
\begin{equation}
 A_i^{\pm} = \left \{ z \in (\RR^k)^n \ \ ; \ \log |z_i| = \pm \log
R, \ \ R^{-1} < |z_j| < R \ \text{for all} \ j \neq i \right \}.
\label{eq_340}
\end{equation}
We write $\cD_R$ for the collection of all functions $u: \Omega_R
\rightarrow \RR$, smooth up to the boundary, that satisfy Neumann's
condition:
\begin{equation} \langle (\nabla u)_i, z_i \rangle = 0
\quad \quad \quad \text{for any} \ \ i=1,\ldots,n, \ \ \ z \in A_i^{\pm}. \label{eq_346_}
\end{equation}
Here, $\nabla u = \left( (\nabla u)_1, \ldots, (\nabla u)_n \right)
\in (\RR^k)^n$. Let $G = ( O(k) ) ^n$, where $O(k)$ is the group of all
orthogonal transformations in $\RR^k$. The group $G$ acts on $\RR^m
= (\RR^k)^n$, via
$$ g.(z_1,\ldots,z_n) =
(g_1(z_1),\ldots,g_n(z_n)) $$ for $g = (g_1,\ldots,g_n) \in G =
O(k)^n$ and $z = (z_1,\ldots,z_n) \in (\RR^k)^n$. A
 subset $U \subseteq \RR^m$ is {\it $G$-invariant} if
 $g.z \in U$ for any $z \in U, g \in G$. Suppose $U \subseteq \RR^m$
 is $G$-invariant and $f: U \rightarrow \RR$. We say that $f$ is
 $G$-invariant if
 $$ f(g.z)  = f(z) \quad \quad \quad \text{for} \ \ g \in G, \ z \in
 U. $$
We write $\pi^{-1}(\RR^n_+)$ for the collection of all $z \in
(\RR^k)^n$ with $z_i \neq 0$ for all $i$. Assume that $\psi:
\pi^{-1}(\RR^n_+) \rightarrow \RR$ is a smooth function, and denote
by $\nu$ the measure on $\pi^{-1}(\RR^n_+)$ whose density is
$\exp(-\psi)$. For a smooth function $u: \pi^{-1}(\RR^n_+)
\rightarrow \RR$ write
$$ \triangle^{\nu} u = e^{\psi} div( e^{-\psi} \nabla u) = \triangle u - \langle \nabla \psi, \nabla u
\rangle, $$ where $div$ stands for the usual divergence operator in
$\RR^m$. Integrating by parts, we see that for any $u, f: \Omega_R
\rightarrow \RR$ that are smooth up to the boundary,
$$ \int_{\Omega_R} \langle \nabla u, \nabla f \rangle d
\nu = -\int_{\Omega_R} f \left(\triangle^{\nu} u \right)  d \nu + \int_{\partial_{reg} \Omega_R} f \langle \nabla u, N \rangle e^{-\vphi}, $$
where $N$ is the outer unit normal. In particular, when $f: \Omega_R \rightarrow \RR$ is smooth up to the boundary
and $u \in \cD_R$,
\begin{equation} \int_{\Omega_R} \langle \nabla u, \nabla f \rangle d
\nu = -\int_{\Omega_R} f \left(\triangle^{\nu} u \right)  d \nu.
\label{eq_2114}
\end{equation}
 The well-known Bochner
identity states that for any smooth function $u: \Omega_R
\rightarrow \RR$,
\begin{equation}
 \frac{1}{2} \triangle^{\nu} \left| \nabla u \right|^2 = \left \langle \nabla u, \nabla (\triangle^{\nu} u) \right \rangle +
 \sum_{i=1}^{m} |\nabla \partial^i
 u|^2 + \left \langle (\nabla^2 \psi) \nabla u, \nabla u \right
 \rangle,
\label{eq_330}
 \end{equation}
as may be verified directly.

\begin{lemma} Let $R > 1$ and let $u
\in \cD_R$ be a $G$-invariant function.  Then,
$$ \int_{\Omega_R} \left|\triangle^{\nu} u \right|^2 d \nu = \int_{\Omega_R} \sum_{i=1}^{m} |\nabla \partial^i
 u|^2 d \nu + \int_{\Omega_R} \left \langle (\nabla^2 \psi) \nabla u, \nabla u \right
 \rangle d \nu. $$ \label{lem_322}
\end{lemma}

\emph{Proof:} We integrate the identity (\ref{eq_330}) over
$\Omega_R$. From (\ref{eq_2114}),
$$ \frac{1}{2} \int_{\Omega_R}  \triangle^{\nu} \left| \nabla u \right|^2 d \nu + \int_{\Omega_R} \left|\triangle^{\nu} u \right|^2 d \nu = \int_{\Omega_R} \sum_{i=1}^{m} |\nabla \partial^i
 u|^2 d \nu + \int_{\Omega_R} \left \langle (\nabla^2 \psi) \nabla u, \nabla u \right
 \rangle d \nu, $$
 since $u \in \cD_R$. To conclude the lemma, it suffices to show that
$$
\int_{\Omega_R}  \triangle^{\nu} \left| \nabla u \right|^2 d \nu =
0. $$ This would follow from (\ref{eq_2114}) once we show that
$\left| \nabla u \right|^2 \in \cD_R$. Hence, in order to conclude
the lemma, we need to prove that \begin{equation} \left \langle \left( \nabla
\left| \nabla u \right|^2 \right)_i, z_i \right \rangle = 0 \quad
\quad \quad \text{for any} \ \ i=1,\ldots,n, \ \ \ z \in A_i^{\pm}.
\label{eq_2117}
\end{equation}
 So far we did not apply the $G$-invariance of
$u$. It will play a role in the proof of (\ref{eq_2117}). Fix
$i=1,\ldots,n$. Since $u \in \cD_R$, then according to
(\ref{eq_346_}), for $z \in A_i^{\pm}$,
$$ \langle (\nabla u)_i, z_i \rangle = 0.
$$
However, since $u$ is $G$-invariant, then $(\nabla u)_i$ is always a
vector proportional to $z_i$.  We conclude that
\begin{equation} (\nabla u)_i = 0 \ \ \text{on} \ \ A_i^{\pm}. \label{eq_2214} \end{equation}
We may differentiate (\ref{eq_2214}) in the direction of $\nabla u$,
since $\nabla u$ is tangential to $\partial_{reg} \Omega_R$, and
obtain
\begin{equation}
   \left( (\nabla^2 u) \nabla u \right)_i = 0 \ \ \text{on} \ \ A_i^{\pm}. \label{eq_2247}
\end{equation} Observe that
\begin{equation}  \nabla \left| \nabla u \right|^2  = 2
(\nabla^2 u) \nabla u. \label{eq_2132} \end{equation} From
(\ref{eq_2247}) and (\ref{eq_2132}) we deduce (\ref{eq_2117}).
\hfill $\square$

\begin{lemma} Suppose that $\vphi: \RR^n_+ \rightarrow \RR$ is smooth, and that the function
$$ (x_1,\ldots,x_n) \mapsto \vphi(x_1^k,\ldots,x_n^k) $$
is convex in $\RR^n_+$. For $z \in \pi^{-1}(\RR^n_+)$  denote
$\psi(z) = \vphi(\pi(z))$. Then, for any $G$-invariant function $u:
\RR^m \rightarrow \RR$,
\begin{equation}  \left \langle (\nabla^2 \psi) \nabla u, \nabla u \right
 \rangle  \geq 0 \label{eq_255_} \end{equation}
at any point $z \in \pi^{-1}(\RR^n_+)$ in which $u$ is
differentiable. \label{convexity}
\end{lemma}

\emph{Proof:} Fix a point $z = (z_1,\ldots,z_n) \in (\RR^k)^n$ with
$z_i \neq 0$ for all $i$. Then the function
$$ \RR^n_+ \ni (a_1,\ldots,a_n) \mapsto \psi(a_1 z_1,\ldots, a_n
z_n) \in \RR $$ is convex on $\RR^n_+$, by our assumption. In
particular, $\nabla^2 \psi(z)|_E$ is positive semi-definite, where
$$ E = \left \{  \left(a_1 z_1,\ldots, a_n z_n \right) \ ;  \ a_1,\ldots,a_n \in \RR
\right \} \subset \RR^m $$ is an $n$-dimensional subspace. Since $u$
is $G$-invariant and differentiable at $z$, then $\nabla u(z) \in
E$, and (\ref{eq_255_}) follows. \hfill $\square$

\medskip
Write $\nu_R$ for the restriction of $\nu$ to $\Omega_R$. We will
use the following well-known fact from the theory of strongly
elliptic operators on convex domains:

\begin{lemma} Suppose $R > 1$. Let $f: \Omega_R \rightarrow \RR$ be a $G$-invariant function that is smooth up to the boundary with $\int f d \nu_R = 0$.
Then, there exists a $G$-invariant function $u \in \cD_R$ with $\int
u d \nu_R = 0$ such that
\begin{equation} \triangle^{\nu} u = f \quad \quad \quad \text{in} \ \ \Omega_R. \label{eq_600} \end{equation}
\label{pde}
\end{lemma}

\emph{Proof sketch:} Denote $Q_R = [-1/R, R]^n \subset \RR^n$ and
$g(|z_1|,\ldots,|z_n|) = f(z_1,\ldots,z_n)$ for $z \in \Omega_R$.
Then $g$ is smooth up to the boundary in $Q_R$. Denote by $\eta$
the finite Borel measure on $Q_R$ which is the push-forward of the
measure $\nu_R$ under the map $(z_1,\ldots,z_n) \mapsto
(|z_1|,\ldots,|z_n|)$. Then $\eta$ has a density of the form
$\exp(-\theta)$ on $Q_R$, where $\theta$ is smooth up to the
boundary.
 Furthermore, $\int g d \eta = 0$. The task of solving (\ref{eq_600}) is reduced to the task of finding $u: Q_R \rightarrow \RR$,
smooth up to the boundary with $\int u d \eta = 0$, such that
\begin{equation}  \triangle u = g + \left \langle \nabla u, \nabla \theta \right \rangle, \label{eq_831} \end{equation}
and such that $u$ satisfies Neumann's boundary condition on
$\partial Q_R$. First, with the help of a crude Poincar\'e
inequality and the Riesz representation theorem, we find a weak
solution. That is, we find $u$ in the Sobolev space $H^1(Q_R) = W^{1,2}(Q_R)$
with $\int u d \eta = 0$ such that (\ref{eq_831}) holds in the sense
that
\begin{equation}  \int_{Q_R} \left \langle \nabla u, \nabla h \right \rangle d \eta = -\int_{Q_R} g h d \eta \quad \quad \quad \text{for any} \ \ h \in H^1(Q_R).
\label{eq_853}
\end{equation} See, e.g., Brezis \cite[Chapter 9]{brezis} or Folland \cite[Chapter 7]{folland} for further explanations. Since $\theta$ is smooth up to the boundary, then $u \in H^k$ implies
$\left \langle \nabla u, \nabla \theta \right \rangle \in H^{k-1}$
for any $k \geq 1$. Furthermore, by expanding into Fourier series in
the cube $Q_R$, one sees that $\triangle u \in H^k$ implies $u \in
H^{k+2}$ for any $k \geq 0$. Therefore, for any $k \geq 0$, if $u
\in H^k$ then from (\ref{eq_831}) also $\triangle u \in H^{k-1}$,
and hence $u \in H^{k+1}$. Therefore $u \in H^k$ for all $k$, and
$u$ is smooth up to the boundary in $Q_R$. From (\ref{eq_853}) we
deduce that
$$ \int_{Q_R} h \left( \triangle u - g - \left \langle \nabla u, \nabla \theta \right \rangle \right) d \eta = \int_{\partial Q_R} h \langle \nabla u, N \rangle e^{-\theta} $$
for any function $h$ that is smooth up to the boundary in $Q_R$.
Here, $N$ is the outer unit normal. This implies that (\ref{eq_831})
holds true in the classical sense, and that $u$ satisfies Neumann's
condition at $\partial Q_R$, as required. \hfill $\square$

\begin{lemma} Let $\vphi$ be as in Lemma \ref{convexity}.
Suppose that $\mu$ is a Borel measure on $\RR^{n}_+$ with density
$\exp(-\vphi)$. Then, for any locally Lipschitz function $f \in
L^2(\mu) \cap L^1(\mu)$,
\begin{equation}
Var_{\mu}(f)  \leq \frac{k^2}{k-1} \sum_{i=1}^n \int_{\RR^n_{+}}
x_i^{2} \left| \partial^i f(x) \right|^2 d \mu(x). \label{eq_913}
\end{equation}
Here, $Var_{\mu}(f) = \int (f - E)^2 d \mu$, where $E \in \RR$ is
such that $\int (f -E) d\mu = 0$.
 \label{lem_2253}
\end{lemma}

\emph{Proof:} By a standard approximation argument (e.g., convolve
$f$ with a localized bump function), we may assume that $f$ is
smooth on $\RR^n_+$. Denote $\psi(z) = \vphi(\pi(z))$ for  $z \in
\pi^{-1}(\RR^n_+)$. Let $\nu$ be the measure on $\RR^m$ whose
density is
$$ z \mapsto \omega_{n,k}^{-1} \exp(-\psi(\pi(z))) \quad \quad \quad \quad (z \in
\pi^{-1}(\RR^n_+)) $$ where $\omega_{n,k}$ is as in Lemma
\ref{lem_2311}. Then $\pi$ pushes the measure $\nu$ forward  to the
measure $\mu$, as we learn from Lemma \ref{lem_2311}, and in fact,
\begin{equation}
\nu = \int_{\RR^n_{+}} \sigma_x d \mu(x).
 \label{eq_1007}
\end{equation}
Fix $R > 1$ and denote $g(z) = f(\pi(z))$. The function $g$ is
smooth up to the boundary in $\Omega_R$. Let $E_R \in \RR$ be such
that $\int (g - E_R) d \nu_R = 0$. According to Lemma \ref{pde},
there exists a $G$-invariant function $u \in \cD_R$ with $\int u d
\nu_R = 0$ such that $\triangle^{\nu} u = - (g-E_R)$.
 Lemma \ref{lem_322} and Lemma \ref{convexity} imply that
\begin{equation}
 \int_{\Omega_R} |\triangle^{\nu} u|^2 d \nu \geq \int_{\Omega_R}
\sum_{i=1}^{m} |\nabla \partial^i u|^2 d \nu. \label{eq_2302}
\end{equation}
 We repeat
the duality argument from \cite[Section 2]{ptrf}:
\begin{align} \label{eq_1746}
&   \int (g - E_R)^2 d\nu_R \\ & =  -  \int g \triangle^{\nu} u d
\nu_R  = \sum_{i=1}^m \int \partial^i
 g
\partial^i u d \nu_R \leq \sum_{i=1}^m \| \partial^i g \nonumber
\|_{H^{-1}(\nu_R)} \sqrt{ \int |\nabla \partial^i u|^2 d \nu_R} \\
& \leq \sqrt{\sum_{i=1}^m\| \partial^i g \|_{H^{-1}(\nu_R)}^2}
\sqrt{ \int \sum_{i=1}^m  |\nabla \partial^i u|^2 d \nu_R} \leq
\sqrt{\sum_{i=1}^m\|
\partial^i g \|_{H^{-1}(\nu_R)}^2} \sqrt{\int |\triangle^{\nu} u|^2 d
\nu_R }, \nonumber
\end{align}
where we used (\ref{eq_2302}) in the last inequality. Therefore,
\begin{equation} \int_{\Omega_R} (g - E_R)^2 d \nu_R \leq
\sum_{i=1}^{m} \| \label{eq_1725}
\partial^i g \|_{H^{-1}(\nu_R)}^2   = \sum_{\ell=1}^n
\sum_{j = 1}^k \left \| \frac{\partial g}{\partial z_{\ell}^j}
\right \|_{H^{-1}(\nu_R)}^2. \end{equation} According to Lemma
\ref{lem_1009} and to (\ref{eq_1007}), for any $\ell=1,\ldots,n$ and
$j=1,\ldots,k$,
\begin{equation} \left \| \frac{\partial g}{\partial z_{\ell}^j}
\right \|_{H^{-1}(\nu_R)}^2 \leq \int_{\RR^n_{+}}  \left \|
\frac{\partial g}{\partial z_{\ell}^j} \right
\|_{H^{-1}(\sigma_x)}^2 d \mu(x) \leq \frac{k}{k-1} \int_{\RR^n_{+}}
x_{\ell}^2 \left|
\partial^\ell f (x) \right|^2 d \mu(x),
\label{eq_1013} \end{equation} where the last inequality is the
content of Lemma \ref{lem_2141}. By combining (\ref{eq_1725}) and
(\ref{eq_1013}), and letting $R$ tend to infinity, we obtain
$$ Var_{\mu}(f) = Var_{\nu}(g) \leq
 \frac{k^2}{k-1} \sum_{i=1}^n \int_{\RR^n_{+}} x_i^{2} \left|
\partial^i f(x) \right|^2 d \mu(x). $$
 \hfill $\square$

\medskip \emph{Proof of Theorem \ref{almost_gap}:} Assume first that $\vphi$ is finite and smooth.
All we need in order to deduce (\ref{eq_1015}) from
(\ref{eq_913}) is to remove the assumption that $f \in L^2(\mu)$. To
that end, given a locally Lipschitz $f \in L^1(\mu)$ and $M
> 0$, we consider the truncation $$ f_M = \max \{ \min \{ f, M \},
-M \}.
$$
Then $f_M \in L^2(\mu)$ is locally Lipschitz. The set $E_M = \{ x
\in \RR^n ; |f(x)| = M \}$ is of measure zero for almost every $M >
0$, as $E_M \cap E_{\tilde{M}} = \emptyset$ for $M \neq \tilde{M}$.
 We apply
(\ref{eq_913}) for $f_M$ and let $M$ tend to infinity, and obtain
(\ref{eq_1015}). This completes the proof in the case where $\vphi$
is finite and smooth. For the general case, a standard approximation
argument is needed. One possibility is to observe that it is enough
to prove the theorem where the integrals over $\RR^n_+$ are replaced
by integrals over the cube $$ \left[ R^{-1}, R \right]^n \subset
\RR^n_+,
$$ for any $R
> 1$. On the bounded cube, it is straightforward to approximate
$\exp(-\vphi)$ by a finite, smooth density, such that both the
left-hand side and the right-hand side of (\ref{eq_1015}) are
well-approximated, for a given locally Lipschitz function $f$. This
completes the proof. \hfill $\square$

\begin{remark}{\rm Suppose $k_1,\ldots, k_n \geq 2$ are integers, and that the function $\vphi: \RR^n_+ \rightarrow (-\infty, \infty]$ is such that
$$ (x_1,\ldots,x_n) \mapsto \vphi(x_1^{k_1},\ldots, x_n^{k_n}) $$
is convex on $\RR^n_+$. It is straightforward to adapt the proof of
Theorem \ref{almost_gap} to this case. We obtain a variant of
Theorem \ref{almost_gap}, in which the inequality (\ref{eq_1015}) is
modified as follows: The factor $k^2 / (k-1)$ is inserted into the
sum, and replaced by $k_i^2 / (k_i - 1)$. See Theorem
\ref{thm_aftermath} below.
 \label{rem_1237}
}\end{remark}

\section{Toric K\"ahler Manifolds}
\label{sec3}

This section  provides a
proof of Theorem \ref{thm_2253}. Throughout this section, we assume
that we are given a convex body $K \subset \RR^n$, and a smooth,
convex function $\psi: \RR^n \rightarrow \RR$ with $\nabla
\psi(\RR^n)= K$. Most of the argument generalizes to any open, convex set $K \subset \RR^n$. In particular, the analysis in Section \ref{sec2} for $k=2$ is parallel to
the case where $K$ equals $\RR^n_+$ and $\psi(x) =
\sum_{i=1}^n \exp(x_i)$.

\medskip The  proof of Theorem
\ref{thm_2253} is essentially  an interpretation of the dual Bochner
inequality in a certain toric K\"ahler manifold. We begin with a
quick review of the the basic definitions, see e.g. Tian
\cite[Chapter 1]{tian} for more information. Suppose $X$ is a
complex manifold of complex dimension $n$. The induced almost
complex structure is a certain smooth map  $J: T X \rightarrow T X$,
such that for any $p \in X$ the restriction $J|_{T_p X}$ is a linear
operator onto $T_p X$ with
$$ J^2|_{T_p X} = - I. $$
In fact, in an open set $U \subset \CC^n$ containing the origin,
consider the map $f(z) = \sqrt{-1} \,  z$ defined in a neighborhood
of zero. Its derivative at zero is $J|_{T_0 U}$. One verifies that
this construction of $J$ does not depend on the choice of the chart,
as the transition functions are holomorphic.  A closed $2$-form
$\omega$ on $X$ is {\it K\"ahler} if the bilinear form
$$ g_{\omega}(u, v) = \omega(u, J v) \quad \quad \quad \quad (p \in X, \ \ \ u,v \in T_p X) $$
is a Riemannian metric, which is also $J$-invariant (i.e.,
$g_{\omega}(u,v) = g_{\omega}(J u, J v)$ for any $p \in X$ and $u,v
\in T_p X$). Next, we specialize to the case of toric K\"ahler
manifolds, see also Abreu \cite{abr} and Gromov \cite{gromov}. We
consider the complex torus
$$ \TT_{\CC}^n = \CC^n / (\sqrt{-1} \ZZ^n) = \left \{ x + \sqrt{-1} y \ ; \ x \in
\RR^n, \ y \in \RR^n / \ZZ^n \right \}. $$ (Perhaps it is more
common to say that $(\CC^*)^n$ is the complex torus, where $\CC^* =
\CC \setminus \{ 0 \}$. Note that $\exp(2 \pi z)$ is a
biholomorphism between $\TT_{\CC}^1$ and $\CC^*$). The real torus
$\TT^n = \RR^n / \ZZ^n$ acts on the complex manifold $\TT_{\CC}^n$
via
$$ t. (x + \sqrt{-1} y) = x + \sqrt{-1} (y + t) \quad \quad \quad
\left(t \in \TT^n, x+ \sqrt{-1} y \in \TT_{\CC}^n \right). $$
Functions, vector fields and differential forms on $\RR^n$ have
toric-invariant extensions to $\TT_{\CC}^n$. For instance, we extend
the convex function $\psi$ to $\TT_{\CC}^n$ by
$$ \psi(x + \sqrt{-1} y ) = \psi(x) \quad \quad \quad \quad \text{for} \
x+\sqrt{-1}y \in \TT_{\CC}^n. $$ Then $\psi$ is a $\TT^n$-invariant
function on the complex manifold $\TT_{\CC}^n$. With a slight abuse
of notation, we use the same letter to denote a function on $\RR^n$,
and its toric-invariant extension to $\TT_{\CC}^n$. Consider the
K\"ahler form on $\TT_{\CC}^n$ defined by
$$ \omega_{\psi} = 2 \sqrt{-1} \partial \bar{\partial} \psi = \frac{\sqrt{-1}}{2} \sum_{i,j=1}^n \psi_{ij} dz_i \wedge
d\bar{z}_j. $$ Abbreviating $g_{\psi} = g_{\omega_{\psi}}$, we have
$$ g_{\psi} \left( \frac{\partial}{\partial x_i},
\frac{\partial}{\partial x_j} \right) = g_{\psi} \left(
\frac{\partial}{\partial y_i}, \frac{\partial}{\partial y_j} \right)
= \psi_{ij} \quad \quad \quad (i,j=1,\ldots,n) $$ while $g_{\psi}
\left( \frac{\partial}{\partial x_i}, \frac{\partial}{\partial y_j}
\right) = 0$ for any $i,j$.
 Furthermore,
observe that
$$ \omega_{\psi}^n = \rho_{\psi} Vol_{2n} $$
where $Vol_{2n}$ is the standard volume form on $\TT_{\CC}^n$ and
$\rho_{\psi}(x) = \det \nabla^2 \psi(x)$ for $x \in \RR^n$.
 It is customary to call the map $x + \sqrt{-1} y
\mapsto \nabla \psi(x)$ the associated {\it moment map}, see Abreu
\cite{abr} and Gromov \cite{gromov}.

\medskip Below we review  in great detail some of the standard formulae of Riemannian geometry
in the case of a toric K\"ahler manifold. As much as possible, we
prefer real formulae in real variables. One reason for this is that
the complex notation fits well only with the case $k=2$ in Section
\ref{sec2}. For a smooth function $u: \RR^n \rightarrow \RR$ we
write
$$ \nabla^{\psi} u = \sum_{i,j=1}^n \psi^{i j} u_i
\frac{\partial}{\partial x_j} = \sum_{j=1}^n u^j
\frac{\partial}{\partial x_j} $$ for the Riemannian gradient of $u$,
where we abbreviate $u^j = \sum_{i=1}^n \psi^{ij} u_i$. Next, we
describe the connection $\nabla^{\psi}$ that corresponds to the
Riemannian metric $g_{\psi}$. As is computed, e.g., in Tian
\cite{tian},
$$ \nabla^{\psi}_{\frac{\partial}{\partial y_j}} \frac{\partial}{\partial
x_k} = \frac{1}{2} \sum_{\ell=1}^n \psi_{jk}^{\ell}
\frac{\partial}{\partial y_\ell}, \quad \quad
\nabla^{\psi}_{\frac{\partial}{\partial x_j}}
\frac{\partial}{\partial x_k} = \frac{1}{2} \sum_{\ell=1}^n
\psi_{jk}^{\ell} \frac{\partial}{\partial x_\ell} $$ where
$\psi_{jk}^{\ell} = \sum_{m=1}^n \psi^{\ell m} \psi_{jkm}$. We view
the Hessian $\nabla^{\psi, 2} h$ of a smooth function $h: \RR^n
\rightarrow \RR$ as a linear operator on $T_p X$, specifically,
$$ T_p X \ni U \mapsto \nabla^{\psi}_{U} \nabla^{\psi} h \in T_p X.
$$
 In
coordinates, for a smooth function $h: \RR^n \rightarrow \RR$, $$
\nabla^{\psi, 2} h \left( \frac{\partial}{\partial x_i} \right) =
\sum_{j,k=1}^n \left( \psi^{jk} h_{ik} - \frac{1}{2} \psi_i^{jk} h_k
\right) \frac{\partial}{\partial x_j}, $$
\begin{equation} \nabla^{\psi, 2} h \left( \frac{\partial}{\partial y_i} \right) =  \frac{1}{2} \sum_{j,k=1}^n \psi_i^{jk} h_k \frac{\partial}{\partial
y_j}, \label{eq_522}
\end{equation}
where $\psi_{i}^{jk} = \sum_{\ell,m=1}^n \psi^{\ell j} \psi^{m k}
\psi_{i\ell m}$. It is unfortunate that we have to work with the
real Hessian, and not with the simpler complex Hessian. We denote by
$\triangle^{\psi}$ the Riemmanian Laplacian on $\TT_{\CC}^n$,
corresponding to the Riemmanian metric $g_{\psi}$. Then
$\triangle^{\psi} h$ is the trace of $\nabla^{\psi, 2} h$, and for a
smooth function $h: \RR^n \rightarrow \RR$,
$$ \triangle^{\psi} h = \sum_{i,j=1}^n \psi^{ij} h_{ij}. $$
The Bochner-Weitzenb\"ock formula from Riemannian geometry (e.g.
Petersen \cite[Section 7.3.1]{petersen}) states that for any smooth
function $u: \RR^n \rightarrow \RR$, \begin{equation} \frac{1}{2}
\triangle^{\psi} |\nabla^{\psi} u|^2 = \langle \nabla^{\psi} u,
\nabla^{\psi} (\triangle^{\psi} u) \rangle + |\nabla^{{\psi},2}
u|_{HS}^2 + Ric_{\psi}(\nabla^{\psi} u, \nabla^{\psi} u)
\label{eq_1218}
\end{equation} where $|\nabla^{{\psi},2} u|_{HS}^2$ is the
Hilbert-Schmidt norm of the Hessian, and where $Ric_{\psi}$ is the
Ricci form, which is the bilinear form given by
$$ Ric_{\psi} \left( \frac{\partial}{\partial x_j} , \frac{\partial}{\partial x_k} \right) = -\frac{1}{2} \frac{\partial^2 \log \rho_{\psi}}{\partial x_j \partial x_k} $$
for $j,k=1,\ldots,n$. Note that $Ric_{\psi}(\nabla^{\psi} u,
\nabla^{\psi} u) \geq 0$ when $\rho_{\psi}$ is log-concave.

\begin{definition} Suppose $(M, g)$ is a Riemannian manifold,
$\nabla$ is the standard Levi-Civita connection, and $\nu$ a Borel
measure on $M$. Let $V$ be a vector field on $M$, which is locally
$\nu$-integrable. We set
\begin{equation}
 \| V \|_{H^{-1}(\nu)} = \sup \left \{ \int_M \langle V, \nabla h \rangle d \nu \
; \ \int_M |\nabla^2 h|_{HS}^2 d \nu\leq 1 \right \} \label{eq_934}
\end{equation} where the supremum runs over all smooth functions $h: M \rightarrow \RR$
such that $\langle V, \nabla h \rangle$ is $\nu$-integrable.
\end{definition}

The proof of Lemma \ref{lem_1009} immediately generalizes to
\begin{equation} \nu =
\int_{\Omega} \nu_{\alpha} d \lambda(\alpha) \quad \quad \Rightarrow
\quad \quad \| V \|_{H^{-1}(\nu)}^2 \leq \int_{\Omega} \| V
\|_{H^{-1}(\nu_\alpha)}^2 d \lambda(\alpha). \label{eq_222}
\end{equation}

Next, we use the $\TT^n$-invariance and obtain a lower bound for $|
\nabla^{\psi, 2} u|_{HS}^2$ in terms of the first derivatives of $u$.
Suppose that $u: \RR^n \rightarrow \RR$ is a smooth function. Denote
by $E_p \subset T_p X$ the subspace spanned by
$\frac{\partial}{\partial y_j} \ (j=1,\ldots,n)$. As in any
Riemannian manifold, the operator $\nabla^{\psi, 2} u$ is symmetric
with respect to the Riemmannian metric $g_{\psi}$. Furthermore, from
(\ref{eq_522}) we learn that $E_p$ is an invariant subspace of the
operator $\nabla^{\psi, 2} u$, and the matrix representing the
operator $\nabla^{\psi, 2} u|_{E_p}$ in the basis
$\frac{\partial}{\partial y_k} \ (k=1,\ldots,n)$ is
$$ \left ( \frac{1}{2} \sum_{j=1}^n u^j \psi^{\ell}_{jk} \right
)_{k,\ell=1,\ldots,n}. $$
 Consequently,
 \begin{align} \nonumber
 \left| \nabla^{\psi, 2} u \right|_{HS}^2 & \geq \left| \left( \left. \nabla^{\psi, 2} u \right|_{E_p} \right) \right |_{HS}^2
= Trace \left[ \left( \left. \nabla^{\psi, 2} u \right|_{E_p}
\right)^2 \right] \\ & = \frac{1}{4} \sum_{i,j, m,p=1}^n u^i u^j
\psi^p_{jm} \psi^m_{i p}. \label{eq_528}
\end{align}
For $x \in \RR^n$ we denote by $\sigma_x$ the uniform probability
measure on the real torus $\{ x + \sqrt{-1} y \, ; \, y \in \TT^n
\}$. For a vector field $U = \sum_{i=1}^n U^i
\frac{\partial}{\partial x_i}$ set
$$ \tilde{Q}_{\psi, x}(U) = \sup \left \{  \left( \sum_{j=1}^n \psi_{ij} U^j
V^j \right)^2 \ ; \ \frac{1}{4} \sum_{i,j,k,\ell, m,p=1}^n V^i V^j
\psi^{\ell m} \psi_{jkm} \psi^{k p} \psi_{i \ell p} \leq 1 \right
\}, $$ where the supremum runs over all $V^1,\ldots,V^n \in \RR^n$.
Here, $\psi^{\ell m}, \psi_{jkm}$ etc. are evaluated at $x$. Observe
that $\tilde{Q}_{\psi, x}$ is essentially the same quadratic form as
$Q_{\psi, \nabla \psi(x)}$ mentioned in the Introduction. That is,
if $h = f(\nabla \psi(x))$, then
$$ \tilde{Q}_{\psi, x} \left( \nabla^{\psi} h \right) =
Q_{\psi, \nabla \psi(x)}(\nabla f). $$

\begin{lemma} Let $u: \RR^n \rightarrow \RR$. Then, for any $x \in \RR^n$ in which $u$ is
differentiable,
$$ \| \nabla^{\psi} u \|^2_{H^{-1}(\sigma_x)} \leq \tilde{Q}_{\psi, x}(\nabla^{\psi} u). $$
\label{lem_1024}
\end{lemma}

\emph{Proof:} The vector field $\nabla^{\psi} u$ on $\TT_{\CC}^n$ is
$\TT^n$-invariant. It therefore suffices to restrict our attention
to $\TT^n$-invariant functions $h$ in the definition (\ref{eq_934})
of $\| \nabla^{\psi} u \|_{H^{-1}(\sigma_x)}$ (i.e., if $h$ is not
$\TT^n$-invariant, then average it with respect to the
$\TT^n$-action). Suppose that $h: \RR^n \rightarrow \RR$ is a smooth
function. From (\ref{eq_528}),
$$ \int_{\TT_{\CC}^n} | \nabla^{\psi, 2} h |_{HS}^2 d \sigma_x \geq \frac{1}{4} \sum_{i,j,k,\ell, m,p=1}^n h^i h^j \psi^{\ell m} \psi_{jkm} \psi^{k
p} \psi_{i \ell p} $$ where the functions on the right-hand side are
evaluated at the point $x$. Since
$$ \int_{\TT^{\CC}_n} \langle \nabla^{\psi} u, \nabla^{\psi} h \rangle d
\sigma_x = \sum_{i,j=1}^n \psi_{ij} u^i h^j, $$
 the lemma follows from the definition of the $H^{-1}$ norm. \hfill
 $\square$

\medskip Suppose $\vphi: \RR^n \rightarrow \RR$ is
a smooth function on $\RR^n$, with $\inf \vphi > -\infty$. Consider
the finite Borel measure $\mu$ on $\TT_{\CC}^n$ that is induced by
the volume form $\exp(-\vphi) \omega_{\psi}^n$. That is, $\mu$ is
the measure on $\TT_{\CC}^n$ whose density with respect to the
standard Lebesgue measure on $\TT_{\CC}^n$ is
$$ \exp(-\vphi(x)) \rho_{\psi}(x). $$
Observe that
\begin{equation}
 \mu = \int_{\RR^n} \sigma_x  e^{-\vphi(x)} \rho_{\psi}(x) dx.
 \label{eq_1024}
 \end{equation}
 For a smooth function $u: \RR^n \rightarrow \RR$  denote
\begin{equation}
 \triangle^{\mu} u = \triangle^{\psi} u - \sum_{i,j=1}^n \psi^{ij}
u_i \vphi_j. \label{eq_2015} \end{equation} Integrating by parts, we
see that when $u, h: \RR^n \rightarrow \RR$ are smooth functions,
with at least one of them compactly-supported,
\begin{equation}  \int_{\TT_{\CC}^n} h (\triangle^{\mu} u) d \mu = -\int_{\TT_{\CC}^n} \langle \nabla^{\psi} u,
\nabla^{\psi} h \rangle d \mu. \label{eq_1116} \end{equation} We
assume that the following Bakry-\'Emery-Ricci condition holds true:
\begin{enumerate}
\item[($\star$)] For any $x \in \RR^n$,
the matrix
$$ \left( \vphi_{i \ell} - \frac{1}{2} \sum_{k=1}^n \psi_{i \ell}^k
\vphi_k -\frac{1}{2} \frac{\partial^2 \log \rho_{\psi}}{\partial x_i
\partial x_\ell}  \right)_{i,\ell=1,\ldots,n} $$
is positive semi-definite.
\end{enumerate}
Condition ($\star$) is equivalent to the pointwise inequality,
\begin{equation}  \left \langle (\nabla^{\psi, 2} \vphi) U, U
\right \rangle + Ric_{\psi} (U, U) \geq 0 \label{eq_306}
\end{equation} for any vector field of the form $U = \sum_{i=1}^n U^i
\frac{\partial}{\partial x_i}$. In the terminology of Bakry and
\'Emery \cite{BE}, condition ($\star$) means that the
Bakry-\'Emery-Ricci tensor (also known as $\Gamma_2$ or the ``second
carr\'e du champ'') is positive semi-definite, when restricted to
the subspace spanned by $\frac{\partial}{\partial
x_1},\ldots,\frac{\partial}{\partial x_n}$. The only case that is
relevant for Theorem \ref{thm_2253}, is when $\rho_{\psi}$ is
log-concave and $\vphi \equiv 1$. Condition ($\star$) clearly holds
true in this case. Theorem \ref{almost_gap} is related to the case
where $\psi(x) = \sum_{i=1}^n e^{x_i}$, and condition ($\star$)
amounts to the convexity of the function $ \vphi(2 \log x_1,\ldots,2
\log x_n)$ in the interior of $\RR^n_+$.

\medskip As explained in the Introduction, we have to impose certain restrictions
on the behavior of $\psi$ and $\vphi$ at infinity. We say that the
pair of functions $(\psi, \vphi)$ is {\it regular at infinity} if
there exists a linear space $X$ of smooth functions $u: \RR^n
\rightarrow \RR$ which has the following properties:
\begin{enumerate}
\item[(a)] For any $u,h \in X$ we have that $h \triangle^{\mu} u,
\left \langle \nabla^{\psi} u, \nabla^{\psi} h \right \rangle \in
L^1(\mu)$, and the the identity (\ref{eq_1116}) holds true. The same
holds also when $u \in X$, and $h: \RR^n \rightarrow \RR$ is such
that $h(\nabla \psi^*(x))$ is a Lipschitz function on $K$.
\item[(b)] The constant functions belong to $X$. If $u \in X$, then also $\triangle^{\mu} u,
|\nabla^{\psi} u|^2 \in X$.
\item[(c)] Denote by $\cH \subset L^2(\mu)$ the subspace of all functions $f:
\RR^n \rightarrow \RR$ with $\int f d \mu =0$. Then the space
$$ \left \{ \triangle^{\mu} u \ ; \ u \in X \right \} $$ is dense in
$\cH$ in the topology of $L^2(\mu)$.
\end{enumerate}
We say that $\psi$ is regular at infinity if $(\psi, 1)$ is regular
at infinity. Observe that  the space of compactly-supported, smooth
functions might not satisfy (c), as there might exist non-constant,
smooth functions $f \in L^2(\mu)$ with $\triangle^{\mu} f \equiv 0$.
The space $X$ is supposed to capture a sort of ``Neumann's condition
at infinity''. A thorough
 investigation of regularity at infinity is beyond the scope of the present
 paper, which focuses on the Bochner method combined with additional
 symmetries in higher dimension.

\begin{remark} {\rm Suppose that the Riemannian manifold $(\TT_{\CC}^n, g_{\psi})$
admits a smooth compactification. That is, assume that
$(\TT_{\CC}^n, g_{\psi})$ embeds in a compact, smooth Riemannian
manifold $(M, g)$ as a dense subset of full measure,  that the
moment map $\nabla \psi$ extends to a smooth function on the entire
$M$, and that the $\TT^n$-action on $(\TT_{\CC}^n, g_{\psi})$
extends to a $\TT^n$-action on $(M,g)$. In this case, $\psi$ is
regular at infinity: We may define $X$ to be the restriction to
$\TT_{\CC}^n$ of all $\TT^n$-invariant, smooth functions on the
compact Riemannian manifold $M$. Indeed, condition (b) then holds
trivially. As for condition (a), observe that $h$ extends to a
Lipschitz function on $M$ as it is the composition of the Lipschitz
maps $h(\nabla \psi^*)$ and $\nabla \psi$, hence integrations by
parts of $h$ against $\triangle^{\psi} u$ may be carried out in $M$.
We conclude that condition (a) holds true since $\TT_{\CC}^n$ is of
full measure in $M$, and the integrals in (\ref{eq_1116}) are
equivalent to integrals over the entire $M$. Condition (c) follows
from the standard theory of elliptic partial differential equations
on a compact, connected, smooth Riemannian manifold.}
\label{rem_307}
\end{remark}

\begin{remark} {\rm Another relevant type of compactification is related to the
so-called orbifolds or $V$-manifolds, which are smooth manifolds
except for some  rather tame singularities. We refer the reader,
e.g., to Chiang \cite{Ch} for Harmonic analysis on Riemannian
orbifolds. In particular, there is a notion of a smooth function on
the entire orbifold, and the Laplace equation may be solved with
smooth functions on compact orbifolds. We conclude that the function
$\psi$ is regular at infinity whenever $(\TT_{\CC}^n, g_{\psi})$
embeds in a compact Riemannian orbifold as a dense subset of full
measure, such that $\nabla \psi$ and the toric action extend
smoothly to the entire Riemannian orbifold. In the case of $K$ being
a rational, simple polytope, all functions $\psi$ admitting such
embedding were characterized by Abreu \cite{abr_orbi}. He gave a
clear criterion in terms of $\psi^*$, which seems to hold in most
cases of interest. Since rational, simple polytopes are dense among
convex bodies, one is tempted to conjecture that Abreu's mild
condition for regularity at infinity may be generalized to the class
of all convex bodies.
 }\end{remark}

The following lemma is a well-known Bochner-type integration by
parts formula. For completeness, we include its proof.

\begin{lemma} Assume that ($\star$) holds true, and that $(\psi,
\vphi)$ is regular at infinity. Then for any $u \in X$,
$$ \int_{\TT_{\CC}^n} |\triangle^{\mu} u|^2 d \mu \geq
\int_{\TT_{\CC}^n} |\nabla^{{\psi},2} u|_{HS}^2 d \mu.
$$
\label{lem_0026}
\end{lemma}

\emph{Proof:} From  (\ref{eq_1218}) and (\ref{eq_2015}) we obtain
the  identity
\begin{align} \label{bochner_before} & \frac{1}{2} \triangle^{\mu} |\nabla^{\psi} u|^2
\\ & =  \langle \nabla^{\psi} u, \nabla^{\psi} (\triangle^{\mu} u)
\rangle + |\nabla^{{\psi},2} u|_{HS}^2 + Ric_{\psi}(\nabla^{\psi} u,
\nabla^{\psi} u) + \left \langle \left( \nabla^{\psi, 2}  \vphi \right) \nabla^{\psi} u,
\nabla^{\psi} u \right   \rangle. \nonumber
\end{align}
From our assumption ($\star$), \begin{equation}
\label{bochner_after} \frac{1}{2} \triangle^{\mu} |\nabla^{\psi}
u|^2 \geq \langle \nabla^{\psi} u, \nabla^{\psi} (\triangle^{\mu} u)
\rangle + |\nabla^{{\psi},2} u|_{HS}^2. \end{equation} Integrating
the above inequality over $\TT_{\CC}^n$, we obtain
$$ 0 \geq -\int_{\TT_{\CC}^n} |\triangle^{\psi} u|^2 d \mu +
\int_{\TT_{\CC}^n} |\nabla^{{\psi},2} u|_{HS}^2 d \mu, $$ since
$\int_{\TT_{\CC}^n} (\triangle^{\psi} h) d \mu = 0$ for any $h \in
X$. \hfill $\square$

\medskip
Theorem \ref{thm_2253} is the case $\vphi \equiv 1$ of the next
 proposition.

\begin{proposition} Let $K \subset \RR^n$ be a convex body. Suppose that
$\psi, \vphi: \RR^n \rightarrow \RR$ are  smooth  functions, such
that $\psi$ is convex with $\det \nabla^2 \psi(x) > 0$ for any $x
\in \RR^n$, and such that $\inf \vphi > -\infty$. Assume that
$\nabla \psi(\RR^n) = K$, that condition ($\star$) above holds true,
and that $(\psi, \vphi)$ is regular at infinity. Let $\mu$ be the measure (\ref{eq_1024})
and denote by $\nu$ the
finite Borel measure on $K$ which is the push-forward of $\mu$ under
$\nabla \psi$. Then, for any Lipschitz function $f: K \rightarrow
\RR$,
\begin{equation} \int_{K} f d \nu = 0 \quad \quad \Rightarrow \quad \quad \int_{K} f^2 d \nu \leq
\int_{K} Q_{\psi,x}(\nabla f) d \nu. \label{eq_223}
\end{equation}
 \label{thm_2253_}
\end{proposition}

\medskip \emph{Proof:} We denote $h(x) = f(\nabla \psi(x))$. Let $u \in X$. With the help
of Lemma \ref{lem_0026}, the duality argument (\ref{eq_1746}) is
replaced by
\begin{align} \label{eq_255}
& -\int_{\TT_{\CC}^n} h \left( \triangle^{\mu} u \right) d \mu   =
\int_{\TT_{\CC}^n} \langle \nabla^{\psi}
 h, \nabla^{\psi} u \rangle
d \mu  \\ & \leq \| \nabla^{\psi} h \|_{H^{-1}(\mu)} \sqrt{
\int_{\TT_{\CC}^n} \left| \nabla^{\psi,2} u \right|_{HS}^2 d \mu }
\leq \| \nabla^{\psi} h \|_{H^{-1}(\mu)} \sqrt{\int_{\TT_{\CC}^n}
|\triangle^{\mu} u|^2 d \mu }. \nonumber
\end{align}
Since $f$ is bounded, then also is $h$ is bounded, hence $h \in
L^2(\mu)$ with $$ \int_{\TT_{\CC}^n} h d \mu = \int_K f d \nu = 0.
$$ Consequently, there exists $u_k \in X$ for $k=1,2,\ldots$ such that $\triangle^{\mu} u_k
\rightarrow -h$ when $k \rightarrow \infty$, in the topology of
$L^2(\mu)$. From (\ref{eq_255}),
$$ \int_K f^2 d \nu = \int_{\TT_{\CC}^n} h^2 d \mu \leq \| \nabla^{\psi}
h
\|_{H^{-1}(\mu)}^2. $$ Combine the latter inequality with
(\ref{eq_222}), (\ref{eq_1024}) and Lemma \ref{lem_1024}, and obtain
\begin{align*} \label{eq_253} \int_{K} f^2 d \nu
 & \leq \| \nabla^{\psi} h
\|_{H^{-1}(\mu)}^2 \leq \int_{\RR^n} \|
\nabla^{\psi} h \|_{H^{-1}(\sigma_x)}^2 e^{-\vphi(x)} \rho_{\psi}(x) dx \\
& \leq \int_{\RR^n} \tilde{Q}_{\psi, x} \left(  \nabla^{\psi} h
\right) e^{-\vphi(x)} \rho_{\psi}(x) dx = \int_{K} Q_{\psi, x}
\left( \nabla f \right) d \nu(x). \nonumber
\end{align*}
\hfill $\square$

\begin{remark}{\rm In principle, one may formulate and prove Theorem
\ref{thm_2253} in terms of $\psi^*$, rather than going back and
forth between  $\psi$ and $\psi^*$, or between $\RR^n$ and $K$. The
reason for preferring $\psi$, is that for $n > 1$, the condition
that $\psi$ induces a log-concave transportation for $K$ appears
simpler than the corresponding condition for $\psi^*$. On the other hand, for a convex function $\psi$ in one variable,
$\log \psi^{\prime \prime}$ is concave if and only if $1 /
(\psi^*)^{\prime \prime}$ is concave. }\end{remark}

\begin{remark}{\rm When $(X, \mu, d)$ is a metric measure space
and $T: X \rightarrow Y$ is a locally Lipschitz map, we may
trivially  transfer any Poincar\'e type inequality on $X$ to a
Poincar\'e type inequality on $Y$. An example is given in Corollary
\ref{cor_456} below, where a Poincar\'e type inequality for the
simplex is deduced from the standard Poincar\'e inequality on $\CC
\PP^n$. Similarly, when $\rho_{\psi} = \exp(-|x|^2 / 2)$, we may, in
principle, transfer the standard Poincar\'e inequality of the
gaussian measure to an inequality on $K$. The approach that we
promote in this paper, of using ``dual Bochner in a higher dimension
with extra symmetries'', is different, and it seems to be applicable
to situations in which the former method fails. Note that we do not
assume any Poincar\'e-type inequality for the log-concave density
$\rho_{\psi}$. } \label{rem_950} \end{remark}

\section{An Example: The Simplex}
\label{sec4}

In order to demonstrate the potential of our paradigm, we present in
this section the Poincar\'e-type inequalities that follow from
Theorem \ref{thm_2253} in the particular case of the simplex. We also discuss
the inequalities that follow via the direct method outlined in
Remark \ref{rem_950}. Our first goal is to apply Theorem
\ref{thm_2253} in the setting
 where $K \subset \RR^n$ is the open simplex whose
vertices are $0, e_1,\ldots, e_n \in\RR^n$. Here, $e_1,\ldots,e_n$
are the standard unit vectors in $\RR^n$. Note that this simplex is
not regular; Later, we will translate the results to the regular
simplex. Consider the smooth, convex function,
$$ \psi(x_1,\ldots,x_n) = \log \left( 1 + e^{x_1} + \ldots + e^{x_n}
\right) \quad \quad \quad (x \in \RR^n). $$ Note that
\begin{equation}
 \nabla \psi(x) = \frac{\left( e^{x_1}, \ldots, e^{x_n} \right) }{1
+ e^{x_1} + \ldots + e^{x_n}}. \label{eq_1021} \end{equation} It is
straightforward to verify from (\ref{eq_1021}) that $$ \nabla
\psi(\RR^n) = K.
$$ Our choice of $\psi$ is motivated by the fact that the K\"ahler
manifold $(\TT_{\CC}^n, \omega_\psi)$ is isometric, up to a
normalization, to a dense open subset of full measure of the complex
projective space $\CC \PP^n$ with the Fubini-Study metric, see e.g.,
the first pages of Tian \cite{tian} or Cannes da Silva \cite{c} for
more information. For instance, the Riemannian manifold
$(\TT_{\CC}^1, g_{\psi})$ is precisely the two-dimensional sphere of
radius one, without the north and the south poles. The moment map
$\nabla \psi$ and the toric action may be extended smoothly to $\CC
\PP^n$, and in view of Remark \ref{rem_307}, we deduce that the
function $\psi$ is regular at infinity. We continue by computing the
second derivatives,
$$ \nabla^2 \psi(x) = \left ( \frac{e^{x_i} \delta_{ij} }{1 + e^{x_1} + \ldots +
e^{x_n}} - \frac{e^{x_i + x_j} }{\left( 1 + e^{x_1} + \ldots +
e^{x_n} \right)^2} \right)_{i,j=1,\ldots,n}. $$ Here, $\delta_{ij}$
is Kronecker's delta.
\begin{lemma} \begin{itemize}
\item[(a)] The function
$$ x \mapsto \det \nabla^2 \psi(x) $$
is log-concave in $\RR^n$.
\item[(b)] The inverse hessian matrix is
$$ \psi^{ij}(x)  = \left(1 + \sum_{j=1}^n e^{x_j} \right) \left[ 1 +
\delta_{ij} e^{-x_i} \right]. $$
\end{itemize}
\end{lemma}

\emph{Proof:} Denote
$$ v = \frac{\left( e^{x_1},\ldots, e^{x_n} \right)}{1 + e^{x_1} + \ldots +
e^{x_n}} \in \RR^n. $$ We write $$ \nabla^2 \psi(x) =  A - B , $$
where $A$ is a diagonal matrix with $v_i$ at the $i^{th}$ diagonal
entry, and $B = ( v_i v_j )_{i,j=1,\ldots,n} $. The determinant of a
rank-one perturbation has a simple formula:
$$ \det \nabla^2 \psi(x) = \det(A - B) = \det(A) \left[ 1 - \langle A^{-1} v, v \rangle \right]. $$
This boils down to \begin{equation}  \det \nabla^2 \psi(x) =  \exp
\left( -(n+1) \psi(x) + \sum_{j=1}^n x_j \right), \label{eq_953}
\end{equation}
 which is
log-concave as $\psi$ is convex. It remains to prove (b). According
to the Sherman-Morisson formula for the inverse of a rank-one
perturbation,
$$ \left( \nabla^2 \psi(x) \right)^{-1} = (A - B)^{-1} = A^{-1} + \frac{A^{-1} B A^{-1}}{1 - \langle A^{-1} v, v \rangle},   $$
as may be verified directly. Equivalently,
$$ \psi^{ij}  = \left(1 + \sum_{j=1}^n e^{x_j} \right) \left[ 1 + \delta_{ij}
e^{-x_i} \right]. $$ \hfill $\square$

\medskip Thus $\psi$ induces a log-concave transportation to $K$.
Note that  $2 Ric_{\psi} = (n+1)  g_{\psi}$, as follows from
(\ref{eq_953}). In particular, we have a very good uniform lower
bound for the Ricci curvature, which implies a rather strong
Poincar\'e inequality on $\CC \PP^n$ -- even a log-Sobolev
inequality -- according to Bakry and \'Emery \cite{BE}.
Consequently, the simple, direct method of Remark \ref{rem_950} has
the potential to produce interesting inequalities in the case of the
simplex. Still, first we would like to test the applicability of
Theorem \ref{thm_2253} here, and to that end, we will write down
explicit expressions for the formidable quadratic form $Q_{\psi,
x}$. We compute
that \begin{align*} \psi_{ijk} & = 2 e^{x_i + x_j + x_k - 3 \psi} + e^{x_i - \psi} \delta_{ij} \delta_{jk} \\
& - \left[ e^{x_j + x_k - 2 \psi} \delta_{ij}  + e^{x_i + x_j -
2\psi} \delta_{ik}  +  e^{x_i + x_k - 2 \psi} \delta_{jk}  \right].
\end{align*}
Therefore,  $$ \psi^{\ell}_{jk}  = \sum_{i=1}^n \psi^{i \ell}
\psi_{i j k} =  \delta_{jk} \delta_{j \ell} - \delta_{j \ell} e^{x_k
- \psi}  - \delta_{k \ell} e^{x_j - \psi}
$$
and, for any fixed $i,j = 1,\ldots,n$,
$$
\sum_{k,\ell=1}^n \psi^{\ell}_{jk} \psi^{k}_{i \ell}  = (n+3) e^{x_i
+ x_j -2 \psi} - e^{x_i - \psi} - e^{x_j - \psi} + \delta_{ij} (1 -
2 e^{x_i - \psi} ). $$ Consequently,
\begin{align*}
 Q^*_{\psi, \nabla \psi(x)}(V) & = \sum_{i,j=1}^n V^i V^j \left[ (n+3) e^{x_i + x_j
-2 \psi} - e^{x_i - \psi} - e^{x_j - \psi} + \delta_{ij} (1 - 2
e^{x_i - \psi} ) \right] \\ & = \sum_{i,j,k=1}^n \psi_{ij} a_{k}^i
V^{k} V^j, \end{align*} where, for $i,k=1,\ldots,n$,
$$ a_{k}^i = e^{x_k} \left(1 - e^{-x_i} \right) + \delta_{ik} \left( e^{\psi - x_i} - 2 \right). $$
We are not confused by the minus signs, and we remember that
$Q^*_{\psi, \nabla \psi(x)}$ must be a positive semi-definite
quadratic form on $\RR^n$. Consider for a moment the scalar product
$$ (U, V) = \sum_{i,j=1}^n \psi_{ij} U^i V^j \quad \quad \quad \quad (U, V \in \RR^n) $$
and the linear operator
$$ A(U)  = \left( \sum_{k=1}^n  a_k^i U^k \right)_{i=1,\ldots,n} \in
\RR^n  \quad \quad \text{for} \ \ U = (U^1,\ldots,U^n) \in \RR^n.$$
Then $A$ is symmetric with respect to the scalar product $( \cdot,
\cdot)$, and $ Q^*_{\psi, \nabla \psi(x)}(V) = \left( A (V), V
\right)$ for $V \in \RR^n$. Observe  that  $$ Q_{\psi, \nabla
\psi(x)}(U) = \sup \left \{ 4 (U, V)^2 \ ; \ V \in \RR^n, \
Q^*_{\psi, \nabla \psi(x)}(V) \leq 1 \right \} = 4 \left( A^{-1}(U),
U \right).
$$ Denote $B = A^{-1} = \left( b^{i}_{j} \right)_{i,j=1,\ldots,n}$.
In order to compute the $b^{i}_{j}$'s, we apply the Sherman-Morisson
formula again, and obtain the expression
$$ b_j^i = \frac{\delta_{ij}}{\psi_j^{-1} - 2} - \frac{\psi_j}{\psi_j^{-1}-2} \cdot \frac{e^{\psi} - \psi_i^{-1}}{\psi_i^{-1} - 2} \left( 1 + \sum_{k=1}^n \frac{e^{\psi}
\psi_k - 1}{\psi_k^{-1} - 2} \right)^{-1}. $$ Therefore,
$$
 \sum_{\ell=1}^n \psi_{i \ell} b^{\ell}_j  = \frac{\psi_i^2}{1 - 2
\psi_i} \delta_{ij} + \frac{\psi_i^2}{1 - 2 \psi_i} \cdot
\frac{\psi_j^2}{1 - 2 \psi_j} \cdot \frac{2 - e^{\psi} }{ 1 +
\sum_{k=1}^n \left[ ( e^{\psi} \psi_k - 1 ) / (\psi_k^{-1} - 2)
\right]}.
$$
Finally, recalling that $\psi_{i}, \exp(\psi)$ are to be evaluated
at the point $\nabla \psi^*(x) = (\nabla \psi)^{-1} x$, we obtain
the positive semi-definite quadratic form
\begin{equation}
\frac{1}{4} Q_{\psi, x}(U)  =  \sum_{i=1}^n \frac{x_i^2 |U^i|^2}{1 -
2 x_i} - \left( \sum_{i= 1}^n \frac{x_i^2 U^i}{1 - 2 x_i} \right)^2
\left( \sum_{k=0}^n \frac{x_k^2}{1 - 2 x_k} \right)^{-1}
\label{eq_432}
\end{equation} where we define $x_0 = 1 - \sum_{j=1}^n x_j$. In
conclusion, so far we have obtained the following:

\begin{corollary} Let $K \subset \RR^n$ be the simplex which is  the convex hull
of $0,e_1,\ldots,e_n$, where $e_1,\ldots,e_n$ are the standard unit
vectors in $\RR^n$. Then for any Lipschitz function $f: K
\rightarrow \RR$ with $\int_K f = 0$,
$$ \int_K f^2(x) dx \leq 4 \int_K \left[ \sum_{i=1}^n \frac{x_i^2 \left| \partial^i f \right|^2}{1 - 2 x_i} - \left(
\sum_{k=0}^n \frac{x_k^2}{1 - 2 x_k} \right)^{-1} \left( \sum_{i=
1}^n \frac{x_i^2  \partial^i f }{1 - 2 x_i} \right)^2 \right] dx
$$
where $x_0 = 1 - \sum_{k=1}^n x_k$. \label{cor_1429}
\end{corollary}

Next, observe that Corollary \ref{cor_510} applies for the uniform
measure on the simplex $K$, with   $\ell=2$. We are
unaware of any advantage of Corollary \ref{cor_1429} over the
inequality that follows from Corollary \ref{cor_510} in this case.
Yet, the importance of Corollary \ref{cor_1429} to us is that it
perhaps demonstrates that the very general Theorem \ref{thm_2253} is
not entirely inapplicable. We continue by translating our results to
the regular simplex.

\medskip Recall that $\RR^{n+1}_+$ is the orthant of all $x \in \RR^{n+1}$
with positive coordinates. Consider the  $n$-dimensional regular
simplex
\begin{equation}  \triangle^n = \left \{ (x_0,\ldots, x_n) \in \RR^{n+1}_+ \  ; \ \sum_{j=0}^n x_j = 1
\right \}. \label{eq_528_} \end{equation} Observe that the
projection
$$ (x_0,\ldots,x_{n}) \mapsto (x_1,\ldots,x_n) $$
is a measure preserving one-to-one correspondence between
$\triangle^n$ and $K$. Let $p \in \triangle^n$, and suppose that $f:
\triangle^n \rightarrow \RR$ is differentiable at $p$. For indices
$i, j = 0,\ldots,n$ we set
$$ E^{ij} f(p) = \left( \frac{\partial}{\partial x_i} - \frac{\partial}{\partial x_j} \right) f (p). $$
Observe that $E_{ij} f(p)$ is well-defined, since the vector field
$\partial / \partial x_i - \partial / \partial x_j$ belongs to the
tangent space $T_p \triangle^n$ for any $p \in \triangle^n$.

\begin{theorem} Let $\triangle^n$ be the simplex (\ref{eq_528_}).
 Then for any Lipschitz function
$f: \triangle^n \rightarrow \RR$ with $\int_{\triangle^n} f = 0$,
$$ \int_{\triangle^n} f^2(x) dx \leq 4 \int_{\triangle^n} \left(
\sum_{k=0}^n \frac{x_k^2}{1 - 2 x_k} \right)^{-1} \sum_{i \neq j}
\frac{x_i^2 x_j^2}{(1 - 2 x_i) (1 - 2 x_j)} \left| E^{ij} f
\right|^2 d  x.
$$
Here, the sum runs over the $n(n+1)/2$ distinct pairs of indices
$i,j \in \{0,\ldots,n\}$. \label{thm_540}
\end{theorem}

\emph{Proof:} For $(x_0,\ldots,x_n) \in \triangle^n$ denote
$$ g(x_1,\ldots,x_n) = f(x_0,\ldots,x_n). $$
Then $g: K \rightarrow \RR$ is a Lipschitz function. We compute that
$$ Q_{\psi, x}(\nabla g(x_1,\ldots,x_n)) = 4 \left( \sum_{k=0}^n \frac{x_k^2}{1 - 2 x_k} \right)^{-1} \sum_{i \neq j}
\frac{x_i^2 x_j^2}{(1 - 2 x_i) (1 - 2 x_j)} \left| E^{ij} f
\right|^2 $$ where $Q_{\psi, x}$ is given by (\ref{eq_432}). The
theorem thus follows from Corollary \ref{cor_1429}. \hfill $\square$

\medskip
We would like to compare Theorem \ref{thm_540} with the push-forward
of the usual Poincar\'e inequality on $\CC \PP^n$ via the moment
map. Recall that  $S^{2n+1}(R) = \{ z \in \CC^{n+1} ; \sum_{i=0}^n
|z_i|^2 = R^2 \}$ is the sphere of radius $R$ in $\CC^{n+1}$,
equipped with the induced Riemannian metric. Recall that the
Riemannian manifold $(\TT_{\CC}^n, g_{\psi})$ is embedded in $\CC
\PP^n$ equipped with the Fubini-Study metric, up to some
normalization. In fact, with respect to the normalization dictated
by $\psi$, we may view the complex projective space $\CC \PP^n$ as a
quotient of the sphere $S^{2n+1}(2) \subset \CC^{n+1}$ by a circle
action. If we extend the map $\nabla \psi$ from $\TT_{\CC}^n$ to
$\CC \PP^n$ by continuity, and then lift it to a circle-invariant
function on $S^{2n+1}(2)$, then we obtain the function
$$ S^{2n+1}(2) \ni (z_0,\ldots,z_n) \mapsto \left(
\frac{|z_1|^2}{4} ,\ldots, \frac{|z_n|^2}{4} \right) \in K.
$$ The manifold $\CC \PP^n$ inherits the Poincar\'e inequality for
even functions on the sphere $S^{2n+1}(2)$ (see, e.g., M\"uller
\cite{muller} for the inequality on the sphere). Consequently, the
standard Poincar\'e inequality on $\CC \PP^n$ is the bound
\begin{equation} \int_{\RR^n} u(x) \rho_{\psi}(x) dx = 0 \quad  \Rightarrow
\quad \int_{\RR^n} u^2(x) \rho_{\psi}(x) dx \leq \frac{1}{n+1}
\int_{\RR^n} |\nabla^{\psi} u(x)|^2 \rho_{\psi}(x) dx,
\label{eq_449}
\end{equation} valid for any function $u: \RR^n \rightarrow
\RR$ for which  $x \mapsto u(\nabla \psi^*(x))$ is Lipschitz. (One
way to make sure that indeed $n+1$ is the first non-zero eigenvalue
of $-\triangle^{\psi}$, is to verify that equality in (\ref{eq_449})
is attained for the eigenfunction $u = \psi_1 - 1 / (n+1)$.)
Translating (\ref{eq_449}) to the simplex $K \subset \RR^n$ via the
moment map $\nabla \psi$, we obtain in a straightforward manner:

\begin{corollary} Let $K \subset \RR^n$ be the simplex which is  the convex hull
of $0,e_1,\ldots,e_n$, where $e_1,\ldots,e_n$ are the standard unit
vectors in $\RR^n$. Then for any Lipschitz function $f: K
\rightarrow \RR$ with $\int_K f = 0$,
$$
 \int_K f^2(x) dx \leq \frac{1}{n+1} \int_K \left[ \sum_{i=1}^n x_i \left| \partial^i f \right|^2 - \left( \sum_{i=
1}^n x_i \partial^i f \right)^2 \right] dx. $$ Equivalently, let
$\triangle^n$ be the simplex (\ref{eq_528_}).
 Then for any Lipschitz function $f: \triangle^n \rightarrow \RR$,
 \begin{equation}
  \int_{\triangle^n} f = 0 \quad \quad \Rightarrow \quad \quad  \int_{\triangle^n} f^2(x) dx \leq \frac{1}{n+1} \int_{\triangle^n} \sum_{i \neq j} x_i  x_j \left| E^{ij} f \right|^2 d  x.
  \label{eq_127}
\end{equation}
Here, the sum runs over the $n(n+1)/2$ distinct pairs of indices
$i,j \in \{0,\ldots,n\}$.  \label{cor_456}
\end{corollary}

\medskip Note that when the dimension $n$ is high, for a random
point $x \in K$ we typically have $x_i \approx \frac{1}{n}$.
Therefore Corollary \ref{cor_456} is not so different from Corollary
\ref{cor_1429}, when the dimension is high, while the latter is less
elegant. Since Corollary \ref{cor_456} has a much shorter proof,
then na\"{i}vely it seems that the general method suggested in
Theorem \ref{thm_2253} is not entirely essential  in the case of the
simplex. In a sense, when proving Corollary \ref{cor_1429} we only
used the fact that $\CC \PP^n$ has a non-negative Ricci form, and we
did not fully exploit the relatively high curvature of $\CC \PP^n$.
The picture is different once we use the freedom to select a
suitable weight function $\exp(-\vphi)$ in Proposition
\ref{thm_2253_}. The following theorem provides a taste of the
Poincar\'e-type inequalities on the simplex that follow from
Proposition \ref{thm_2253_}. Recall the notion of a $p$-convex
function from the Introduction.

\begin{theorem} Let $\triangle^n$ be the simplex (\ref{eq_528_}),
let $q \geq 0$ and let $\vphi: \RR_+^{n+1} \rightarrow \RR$ be a
 $(1/2)$-convex function, smooth up to the boundary in
$\triangle^n$, homogenous of degree $q$. Denote $M = \sup_{x \in \triangle^n} \vphi(x)$, and assume that
\begin{equation} M q \leq n. \label{eq_2246} \end{equation}
(Alternatively, we can assume condition (\ref{eq_331}) below in place of (\ref{eq_2246}).)
 Denote by $\nu$ the finite Borel measure on $\triangle^n
\subset \RR^{n+1}$ whose density with respect to the Lebesgue
measure on $\triangle^n$ is
$$ (x_0, \ldots, x_n) \mapsto \exp \left(-\vphi \left(x_0,\ldots,
x_n \right) \right) \quad \quad \quad \quad (x \in \triangle^n).
$$ Then for any Lipschitz function $f: \triangle^n \rightarrow \RR$
with $\int_{\triangle^n} f d\nu = 0$,
$$ \int_{\triangle^n} f^2(x) d \nu(x) \leq 4 \int_{\triangle^n} \left(
\sum_{k=0}^n \frac{x_k^2}{1 - 2 x_k} \right)^{-1} \sum_{i \neq j}
\frac{x_i^2 x_j^2}{(1 - 2 x_i) (1 - 2 x_j)} \left| E^{ij} f
\right|^2 d  \nu(x).
$$
Here, the sum runs over the $n(n+1)/2$ distinct pairs of indices
$i,j \in \{0,\ldots,n\}$. \label{thm_513}
\end{theorem}

\emph{Proof:} Note that $\vphi$ extends by continuity to the closure
$\overline{\RR^{n+1}_+} \setminus \{0 \}$. Define
$$ f(z_0, \ldots,z_n) = \vphi \left( \frac{|z_0|^2}{4},\ldots,\frac{|z_n|^2}{4} \right) \quad \quad \quad (0 \neq z \in \CC^{n+1}), $$
and observe that $f$ is smooth on $S^{2n+1}(2)$ as $\vphi$ is smooth
up to the boundary in $\triangle^n$. For a point $p \in S^{2n+1}(2)$
we write $E_p \subset T_p (S^{2n+1}(2))$ for the subspace spanned by
the gradients of the functions $|z_0|^2,\ldots,|z_n|^2$ on
$S^{2n+1}(2)$. Arguing as in Lemma \ref{convexity}, we see that
$$ \left \langle (\nabla^2 f)U, U \right \rangle \geq 0 \quad \text{for any} \ \ p \in S^{2n+1}(2), U \in
E_p. $$ From (\ref{eq_2246}),
$$
 \left \langle (\nabla^2 f)U, U \right \rangle + \frac{n - qM}{2}  |U|^2 \geq 0
\quad \text{for any} \ \ p \in S^{2n+1}(2), U \in E_p.
$$
Since $f(p) \leq M$ for any $p \in S^{2n+1}(2)$, then $f$ satisfies
\begin{equation}
 \left \langle (\nabla^2 f)U, U \right \rangle + \frac{n - q f(p)}{2}  |U|^2 \geq 0
\quad \text{for any} \ \ p \in S^{2n+1}(2), U \in E_p.
\label{eq_331}
\end{equation}
The remainder  of the proof is devoted to showing that condition
(\ref{eq_331}) suffices for the application of Proposition
\ref{thm_2253_}. To that end, denote by $\pi: S^{2n+1}(2)
\rightarrow \CC \PP^n$ the quotient map, which associates with any
$z \in S^{2n+1}(2)$ the complex line through the origin that passes
through $z$.  Note that when $p \in S^{2n+1}(2)$ is such that
 $\pi(p)
\in \TT_{\CC}^n$, the subspace $\pi_*(E_p)$ is the linear span of
$\partial / \partial x_1,\ldots, \partial / \partial x_n$. We need to check
 that condition ($\star$) from Section
\ref{sec3} holds true, and  that the pair $$ \left(\psi(x), \vphi\left( \frac{(1, e^{x_1}, \ldots, e^{x_n})}{1 + e^{x_1} + \ldots + e^{x_n}} \right) \right) $$
is regular at infinity. The main observation here is that both
requirements are satisfied when
\begin{equation}
\left \langle \left( \nabla^2_{S^{2n+1}(2)} f \right) U, U \right
\rangle + Ric_{S^{2n+1}(2)}(U, U) \geq 0 \quad \text{for any} \ \ p
\in S^{2n+1}(2), U \in E_p. \label{eq_403}
\end{equation}
Here, $\nabla^2_{S^{2n+1}(2)} f$ stands for the Hessian of $f$ with
respect to the Riemannian metric on $S^{2n+1}(2)$. Indeed, it is
straightforward to verify that the Bakry-\'Emery-Ricci tensor of a
smooth function $g: \CC \PP^n \rightarrow \RR$ is positive
semi-definite on $\pi_*(E_p)$, if and only if the
Bakry-\'Emery-Ricci tensor of  $g \circ \pi: S^{2n+1}(2) \rightarrow
\RR$ is positive semi-definite on $E_p$. Hence (\ref{eq_403})
implies condition ($\star$) from Section \ref{sec3}. The regularity
at infinity is not an issue, as $f \circ \pi^{-1}$ is well-defined
and smooth on the entire $\CC \PP^n$.  Since $Ric_{S^{2n+1}(2)}(U,
U) = n |U|^2/2$ and $f$ is homogenous of degree $2q$, then
(\ref{eq_403}) is equivalent to  (\ref{eq_331}). The theorem is thus
proven. \hfill $\square$

\begin{remark}{\rm
Observe that  the Poincar\'e inequality on $\CC \PP^n$, rendered as
(\ref{eq_449}) above, essentially remains true when we replace the
integrals over the entire $\CC \PP^n$ with integrals over a
geodesically-convex subset of $\CC \PP^n$. This follows from the
Bochner formula, with a slightly weaker constant $2 / (n+1)$ in
place of the factor $1 / (n+1)$ from (\ref{eq_449}). See Escovar
\cite[Theorem 4.3]{escobar} for details and for a better constant.
Consequently, (\ref{eq_127}) remains true, up to a factor of two,
when the integrals over $\triangle^n$ are replaced by integrals over
a compact $K \subset \triangle^n$ for which $\pi^{-1}(K)$ is
geodesically-convex. Here, $\pi: \CC \PP^n \rightarrow
\overline{\triangle^n}$ is the moment map. In the case where $n=1$,
the condition on $K$ means that $K$ is connected, contains one of
the endpoints of the interval $\triangle^1$, and is contained in one
of the halves of the interval $\triangle^1$.}
\end{remark}

\begin{remark}{ \rm Assumption (\ref{eq_2246}) and even the more precise condition (\ref{eq_331})
seem a bit strict. We suspect that this is the fault of the hasty
transition from (\ref{bochner_before}) to (\ref{bochner_after})
above. Perhaps a more subtle analysis, in the spirit of Barthe and
Cordero-Erausquin \cite{bce}, may transform the strict condition
(\ref{eq_2246}) into a parameter incorporated in the resulting
Poincar\'e-type inequality.  }\end{remark}

\begin{remark}{\rm Theorem \ref{thm_540} and its generalization Theorem \ref{thm_513} essentially follow by analyzing
the Fubini-Study metric on $\CC \PP^n$. It seems that there is a
developed theory of ``canonical'' K\"ahler metrics on certain toric
manifolds, and in many cases we even have an everywhere non-negative
Ricci form. Our limited understanding of this theory has  so far
 prevented us from extracting additional meaningful Poincar\'e-type
inequalities. }\end{remark}

\section{From the Orthant to the Full Space}
\label{sec5}

In this section we deduce Theorem \ref{thm_easy} from Theorem
\ref{almost_gap} and from  some essentially known facts.  We
say that an unconditional $\rho:\RR^n \rightarrow \RR$ is increasing
when the restriction $\rho|_{\RR^n_{+}}$ is increasing. We say that
it is decreasing when $x \mapsto -\rho(x)$ is increasing.  The
following lemma begins our analysis of the finite-dimensional space
of functions on $\RR^n$ that are constant on each orthant. Recall
the definition (\ref{eq_1636}) of the $H^{-1}$ norm of a function.

\begin{lemma} Let $R > 0$, and let $\mu$ be the uniform probability measure on the interval
$[-R, R]$. Suppose $f(x) = sgn(x) = x / |x|$ for $x \neq 0$. Then,
\begin{equation}  \| f \|_{H^{-1}(\mu)} \leq \frac{R}{\sqrt{3}} = \sqrt{\int_{\RR} x^2 d \mu(x) }.  \label{eq_938}
\end{equation} \label{lem_1158}
\end{lemma}

\emph{Proof:} Integrating by parts, we see that for any smooth
function $g$,
\begin{align*}
 \frac{1}{2R} \int_{-R}^R f g = \frac{1}{2R} \int_0^R \left[ g(x) - g(-x) \right] dx = \frac{1}{2R} \int_0^R (R - x) \left( g^{\prime}(x) + g^{\prime}(-x) \right) dx \\
 \leq \frac{1}{2R} \sqrt{\int_0^R (R - x)^2 dx \int_0^R \left| g^{\prime}(x) + g^{\prime}(-x) \right|^2 dx} \leq
 \frac{1}{2R} \sqrt{ \frac{2 R^3}{3} \int_{-R}^R \left| g^{\prime}(x) \right|^2 dx },
 \end{align*}
where we used the Cauchy-Schwartz inequality. The bound
(\ref{eq_938}) now follows from the definition (\ref{eq_1636}) of
the $H^{-1}$-norm. \hfill $\square$

\medskip  Suppose $\rho: \RR \rightarrow \RR$ is a probability density that is unconditional (i.e., even) and decreasing.
It is elementary to verify that there exists a probability measure
$\lambda$ on $[0, \infty)$, such that
$$ \rho(x) = \int_0^{\infty} \left(
\frac{1_{[-R, R]}(x) }{2R} \right) d \lambda(R) \quad \quad \quad
(\text{for almost every} \ x \in \RR)
$$
where $1_{[-R, R]}$ is the characteristic function of the interval
$[-R, R]$. From Lemma \ref{lem_1009} and Lemma \ref{lem_1158} we
 conclude that for any probability measure $\mu$ on $\RR$
with an unconditional, decreasing density,
\begin{equation}  \| sgn(x) \|_{H^{-1}(\mu)} \leq  \sqrt{\int_{\RR} x^2 d \mu(x) }.  \label{eq_1537}
\end{equation}
Note that  when $\rho$ is an unconditional, decreasing function on
$\RR^n$,  the restriction of $\rho$ to any line parallel to one of
the axes, is a one-dimensional unconditional, decreasing function.
From (\ref{eq_1537}) and Lemma \ref{lem_1009} we therefore obtain
the following:

\begin{corollary} Suppose $\mu$ is a probability measure on $\RR^n$
with an unconditional, decreasing density. Let $\ell = 1,\ldots,n$,
and suppose that $f: \RR^n \rightarrow \{ -1, 1 \}$ is a measurable
function which does {\it not} depend on the $\ell^{th}$ coordinate.
Set
$$ g(x) = f(x) sgn(x_{\ell}) \quad \quad \quad \quad \text{for} \ x = (x_1,\ldots,x_n) \in \RR^n. $$
Then,
$$ \| g \|_{H^{-1}(\mu)} \leq \sqrt{ \int_{\RR^n} x_{\ell}^2 d \mu(x) }. $$ \label{cor_502}
\end{corollary}

\medskip
Let $G = \{ -1, 1 \}^n \cong (\ZZ / (2 \ZZ))^n$, a commutative group
with $2^n$ elements, where
$$ xy = (x_1 y_1,\ldots,x_n y_n) \quad \quad \text{for} \ \ x,y \in
\{ -1, 1 \}^n. $$ Denote by $\cH$ the space of functions $f: G
\rightarrow \RR$ with $\sum_{x \in G} f(x) = 0$. For $x, y \in G$
and $f \in \cH$ denote $T_x f (y) = f(x y)$.  Suppose that we have
two Hilbertian norms $\| \cdot \|_1$ and $\| \cdot \|_2$ on the
space $\cH$, with the property that
\begin{equation} \| f \|_j = \| T_x f  \|_j \label{eq_452}
\end{equation}
for any $x \in G, f \in \cH$ and $j=1,2$. From elementary
representation theory, the supremum
$$ \sup_{0 \neq f \in \cH} \| f \|_1 / \| f \|_2 $$
must be attained for a non-constant character $f: G \rightarrow
\RR$.

\begin{lemma}  Suppose $\mu$ is a probability measure on $\RR^n$
with an unconditional, decreasing density. Let $\cS \subset
L^2(\mu)$ be the finite-dimensional space spanned by functions $f$
that are constant on orthants. That is, functions $f$ such that
$$ f(x_1,\ldots,x_n) $$
depends only on $sgn(x_1),\ldots,sgn(x_n)$. Then, for any $f \in
\cS$ with $\int f^2 d\mu = 1$ and $\int f d \mu = 0$,
\begin{equation}  \| f \|_{H^{-1}(\mu)}^2 \leq \max_{\ell=1,\ldots,n} \int_{\RR^n} x_\ell^2 d \mu(x). \label{eq_453}
\end{equation} \label{lem_orthant}
\end{lemma}

\emph{Proof:} Denote by $\cH \subset \cS$ the subspace of all
functions $f \in \cS$ with $\int f d \mu = 0$, and consider the
group $G = \{ -1, 1 \}^n \cong (\ZZ / (2 \ZZ))^n$. The linear space
$\cH$ is identified with the space of functions on $G$ that sum to
zero, since each of the $2^n$ orthants is identified with an element
of $G$ in an obvious manner. Furthermore, the $H^{-1}(\mu)$ norm and
the $L^2(\mu)$ norm are both $G$-invariant Hilbertian norms on $\cH$
in the sense of (\ref{eq_452}). It is therefore sufficient to verify
(\ref{eq_453}) for non-constant characters, that is, for functions
$f: \RR^n \rightarrow \RR$ of the form
$$ f(x) = \prod_{j=1}^n sgn(x_j)^{\delta_j} \quad \quad \quad \quad (x \in \RR^n) $$
for some $0 \neq (\delta_1,\ldots,\delta_n) \in \{ 0, 1 \}^n$. Note
that all of these characters are of the form
$$ f(x) = g(x) sgn(x_{\ell}) $$
for some $\ell=1,\ldots,n$ and for some measurable function $g:
\RR^n \rightarrow \{ -1, 1\}$ which does not depend on $x_{\ell}$.
Corollary \ref{cor_502} therefore applies, and implies
(\ref{eq_453}). \hfill $\square$

\medskip
\emph{Proof of Theorem \ref{thm_easy}:} By applying a linear
transformation of the form
$$ \RR^n \ni (x_1,\ldots,x_n) \mapsto (\sqrt{V_1} x_1,\ldots, \sqrt{V_n}
x_n) \in \RR^n
$$
we reduce matters to the case $V_1 = \ldots = V_n = 1$. We will
consider the  norms corresponding to the expressions appearing on
the right-hand side of (\ref{eq_1015}) and of (\ref{eq_0037}). That
is, for a locally Lipschitz function $g \in L^2(\mu)$ set
\begin{align*}
 \| g \|_{P^1(\mu)}^2 & = \int_{\RR^n} \sum_{i=1}^n \frac{k^2}{k-1}
x_i^{2} \left|
\partial^i g(x) \right|^2 d \mu(x), \\
\| g \|_{Q^1(\mu)}^2 & =   \int_{\RR^n} \sum_{i=1}^n \left(
\frac{k^2}{k-1} x_i^{2} + 1 \right) \left|
\partial^i g(x) \right|^2 d \mu(x). \end{align*}
Then
\begin{equation} \| g \|_{Q^1(\mu)}^2 = \| g \|_{P^1(\mu)}^2 + \| g
\|_{H^1(\mu)}^2 \label{eq_1652} \end{equation} where $\| g
\|_{H^1(\mu)}^2 = \int |\nabla g|^2 d\mu$. The dual norms are
defined, for $f \in L^2(\mu)$, via
$$ \| f \|_{P^{-1}(\mu)} = \sup_{\| g
\|_{P^1(\mu)} \neq 0} \frac{\int f g d \mu}{\| g \|_{P^1(\mu)}},
\quad \| f \|_{Q^{-1}(\mu)} = \sup_{\| g \|_{Q^1(\mu)} \neq 0}
\frac{\int f g d \mu}{\| g \|_{Q^1(\mu)}},
$$ where the suprema run over all locally Lipschitz  functions $g \in
L^2(\mu)$. Using a standard duality argument we deduce from
(\ref{eq_1652}) that for any $f_1, f_2 \in L^2(\mu)$,
\begin{equation} \| f_1 + f_2 \|_{Q^{-1}(\mu)}^2 \leq \| f_1
\|_{P^{-1}(\mu)}^2 + \| f_2 \|_{H^{-1}(\mu)}^2 \label{eq_1656}
\end{equation} whenever the right-hand side is finite. In order to
prove (\ref{eq_0037}), it suffices to show that for any $f \in
L^2(\mu)$ with $\int f d\mu = 0$,
\begin{equation}
\| f \|_{Q^{-1}(\mu)} \leq \| f \|_{L^2(\mu)}. \label{eq_1616}
\end{equation}
(Strictly speaking, this will imply (\ref{eq_0037}) only for a
locally Lipschitz $f \in L^2(\mu)$, yet the generalization to a
locally Lipschitz $f \in L^1(\mu)$ is simple, as is explained at the
proof of Theorem \ref{almost_gap} above). For $f: \RR^n
\rightarrow \RR$ and $\delta \in \{-1,1\}^n$ denote
$$ f_{\delta}(x) = f(\delta_1 x_1,\ldots,\delta_n x_n) \quad \quad
\quad \text{for} \ x \in \RR^n_{+}. $$ We write $\cG \subseteq
L^2(\mu)$ for the subspace of all $f \in L^2(\mu)$ which satisfy $$
 \int_{\RR^n_{+}} f_{\delta}
d \mu = 0 \quad \quad \quad \text{for all} \ \delta \in \{ -1, 1
\}^n. $$ Suppose that $g \in L^2(\mu)$ is a locally Lipschitz
function with
\begin{equation}  \| g \|_{P^1(\mu)}^2 = \int_{\RR^n} \sum_{i=1}^n \frac{k^2}{k-1} x_i^{2}
\left|
\partial^i g(x) \right|^2 d \mu(x) \leq 1. \label{eq_1042}
\end{equation}
For $\delta \in \{-1,1\}^n$ let $E_{\delta} \in \RR$ be such that
 $\int_{\RR^n} (g_{\delta} - E_{\delta}) d \mu = 0$.
According to (\ref{eq_1042}) and to Theorem \ref{almost_gap},
$$ \sum_{\delta \in \{ -1,1 \}^n} \int_{\RR^n_{+}} (g_{\delta} - E_{\delta})^2 d
\mu \leq 1. $$ Consequently,  for any $f \in \cG$,
\begin{align*}
 \int_{\RR^n} f g d \mu & = \sum_{\delta \in \{ -1,1 \}^n}
\int_{\RR^n_{+}} f_{\delta} g_{\delta} d \mu  = \sum_{\delta \in \{
-1,1 \}^n} \int_{\RR^n_{+}} f_{\delta} (g_{\delta} - E_{\delta}) d
\mu
\\ & \leq \sqrt{\sum_{\delta \in \{ -1,1 \}^n}  \int_{\RR^n_{+}} f_{\delta}^2 d
\mu } \cdot \sqrt{ \sum_{\delta \in \{ -1,1 \}^n}  \int_{\RR^n_{+}}
(g_{\delta} - E_{\delta})^2 d \mu }\leq \sqrt{ \int_{\RR^n} f^2 d
\mu}. \end{align*} We thus proved that
\begin{equation}
\| f \|_{P^{-1}(\mu)} \leq \| f\|_{L^2(\mu)} \quad \quad \quad \quad
\text{for any} \ f \in \cG. \label{eq_1613}
\end{equation}
Next, observe that $\cG$ is the orthogonal complement to the
subspace $\cS$ from Lemma \ref{lem_orthant}. Fix $f \in L^2(\mu)$
with $\int f d \mu = 0$. Then $f$ may be represented as $ f = g +
s$, where $g \in \cG, s \in \cS$ and $\int s d \mu = 0$. From
(\ref{eq_1656}), (\ref{eq_1613}) and Lemma \ref{lem_orthant},
$$
\| f \|_{Q^{-1}(\mu)}^2 \leq \| g \|^2_{P^{-1}(\mu)} + \| s
\|_{H^{-1}(\mu)}^2 \leq \| g \|_{L^2(\mu)}^2 + \| s \|_{L^2(\mu)}^2
= \| f \|_{L^2(\mu)}^2,
$$
and the desired (\ref{eq_1616}) is proven. The ``Furthermore'' part
of the theorem follows immediately from Theorem \ref{almost_gap}.
\hfill $\square$

\section{A direct approach for the orthant}
\label{aftermath}

In this section we provide another proof of Theorem
\ref{almost_gap}, which does not involve spaces of twice the
dimension. We prove the following slight generalization of Theorem
\ref{almost_gap}, see also Remark \ref{rem_1237}.

 \begin{theorem} Let $n \geq 1$. Let $k_1,\ldots,k_n > 1$ be real numbers, not necessarily integers. Suppose that
 $\mu$ is a Borel measure on $\RR^n_{+}$
with density $\exp(-\vphi)$, where $\vphi: \RR^n_{+} \rightarrow
\RR$ is a smooth function such that
$$ \RR^n_+ \ni (x_1,\ldots,x_n) \mapsto \vphi \left(x_1^{k_1},\ldots,
x_n^{k_n} \right) $$ is a convex function on $\RR^n$.
 Assume  that $f:
\RR^n_+ \rightarrow \RR$ is a $\mu$-integrable, locally Lipschitz
function with $\int f d \mu = 0$. Then,
\begin{equation}
 \int_{\RR^n_{+}} f^2 d \mu \leq
\int_{\RR^n_{+}} \sum_{i=1}^n \frac{k_i^2}{k_i-1} x_i^{2} \left|
\partial^i f(x) \right|^2 d \mu(x). \label{eq_1239} \end{equation}
\label{thm_aftermath}
\end{theorem}

\emph{Proof:}  For $x \in \RR^n_+$  we denote here
$$ \pi(x) = (\pi_1(x),\ldots,\pi_n(x) = (x_1^{k_1},\ldots,x_n^{k_n}). $$
Then $\vphi(\pi(x))$ is a convex function. Set
$$ \psi(x) = \vphi(\pi(x)) - \sum_{i=1}^n (k_i-1) \log x_i \quad \quad \quad \quad (x \in \RR^n_+). $$
Since $\vphi(\pi(x))$ is convex, its Hessian is positive
semi-definite. Therefore,
$$ \left \langle \left( \nabla^2 \psi(x) \right)^{-1} U , U \right \rangle \leq \sum_{i=1}^n \frac{x_i^2}{k_i - 1} |U^i|^2 $$
for any $x \in \RR^n_+$ and $U = (U^1,\ldots, U^n)$. From the
Brascamp-Lieb inequality \cite[Theorem 4.1]{BL}, we conclude that
for any locally Lipschitz function $f: \RR^n_+ \rightarrow \RR$,
\begin{equation}
\int_{\RR^n_+} f e^{-\psi} = 0 \quad \quad \Rightarrow \quad \quad
\int_{\RR^n_+} f^2 e^{-\psi} \leq \int_{\RR^n_+} \sum_{i=1}^n
\frac{x_i^2}{k_i - 1} |\partial^i f(x) |^2 e^{-\psi(x)} dx. \label{eq_1227}
\end{equation}
Equivalently, for any locally Lipschitz  function $f: \RR^n_+
\rightarrow \RR$ with $$ \int_{\RR^n_+} f(x) \left( \prod_{i=1}^n
x_i^{k_i-1} \right) e^{-\vphi(\pi(x))} dx = 0, $$ we have
\begin{equation} \label{eq_1243}
 \int_{\RR^n_+} f^2 \left( \prod_{i=1}^n x_i^{k_i-1} \right)
e^{-\vphi(\pi(x))} dx \leq \int_{\RR^n_+} \sum_{i=1}^n
\frac{x_i^2}{k_i - 1} |\partial^i f |^2 \left( \prod_{i=1}^n
x_i^{k_i-1} \right) e^{-\vphi(\pi(x))} dx. \end{equation} Observe
that $\prod_{i=1}^n k_i x_i^{k_i-1}$ is precisely the Jacobian
determinant of $\pi$. Furthermore, if $f(x) = g(\pi(x))$, then
$$ x_i \partial^i f(x) = k_i \pi_i(x) \partial^i g(\pi(x)). $$
 From (\ref{eq_1243}) we see that for any locally Lipschitz
$f: \RR^n_+ \rightarrow \RR$ with $\int f e^{-\vphi} = 0$,
$$  \int_{\RR^n_+} f^2
e^{-\vphi(x)} dx \leq \int_{\RR^n_+} \sum_{i=1}^n \frac{k_i^2}{k_i -
1} x_i^2 |\partial^i f |^2 e^{-\vphi(x)} dx. $$ \hfill $\square$

\medskip Theorem \ref{thm_aftermath} immediately implies
the corresponding refinements of Corollary \ref{cor_510} and Theorem
\ref{thm_easy}, as described in the Introduction.

\begin{remark}{\rm We currently do not know of any  direct
approach for Theorem \ref{thm_2253} or even for the Poincar\'e
inequalities obtained for the simplex in Section \ref{sec4}. Still,
we cannot escape the feeling that the symmetries we produce by
adding extra dimensions are somewhat artificial. Perhaps we are
overlooking a direct method, that could lead to simpler proofs and
generalizations of the results in this manuscript. }\end{remark}


\end{document}